\documentclass[reqno]{article}

\usepackage{amssymb,amsmath,amscd}

\usepackage[arrow,matrix,curve]{xy}

\newcommand{\gm}{\gamma}
\newcommand{\dl}{\delta}
\newcommand{\sg}{\sigma}
\newcommand{\om}{\omega}

\newcommand{\bR}{\mathbb{R}}

\newcommand{\G}{\,\mathcal{G}}
\newcommand{\cA}{\mathcal{A}}
\newcommand{\cD}{\mathcal{D}}
\newcommand{\cF}{\mathcal{F}}
\newcommand{\cI}{\mathcal{I}}
\newcommand{\cO}{\mathcal{O}}
\newcommand{\cP}{\mathcal{P}}
\newcommand{\cU}{\mathcal{U}}

\newcommand{\Bd}{\,\mathcal{B}d}
\newcommand{\Catg}{\;\mathcal{C}atg}

\newcommand{\Cl}{\,\mathcal{C}l}
\newcommand{\DI}{\,\mathcal{D}\mathcal{I}}
\newcommand{\LC}{\,\mathcal{L}\mathcal{C}}
\newcommand{\ND}{\,\mathcal{N}\mathcal{D}}
\newcommand{\ST}{\,\mathcal{S}\mathcal{T}}

\newcommand{\sbs}{\subset}
\newcommand{\sps}{\supset}
\newcommand{\sbsq}{\subseteq}
\newcommand{\ol}{\overline}
\newcommand{\vnth}{\varnothing}
\newcommand{\ds}{\displaystyle}

\newcommand{\dd}{\mbox{-}\!}
\newcommand{\lla}{\Longleftarrow}
\newcommand{\lra}{\Longrightarrow}
\newcommand{\llra}{\Longleftrightarrow}

\newcommand{\cl}{\,\operatorname{cl}}
\newcommand{\co}{\operatorname{co}}
\newcommand{\nt}{\operatorname{int}}

\newcommand{\R}{\operatorname{R}}
\newcommand{\T}{\operatorname{T}}
\newcommand{\I}{\operatorname{I}}
\newcommand{\II}{\operatorname{II}}

\newcommand{\Fr}{\,\operatorname{Fr}}
\newcommand{\ind}{\,\operatorname{ind}}
\newcommand{\Ind}{\,\operatorname{Ind}}
\newcommand{\lcq}{\,\operatorname{lcq}}
\newcommand{\aind}{\,\operatorname{aind}}
\newcommand{\aInd}{\,\operatorname{aInd}}

\newcommand{\BS}{\,\mbox{\rm BS}}
\newcommand{\TS}{\,\mbox{\rm TS}}
\newcommand{\WS}{\,\mbox{\rm WS}}
\newcommand{\BsS}{\,\mbox{\rm BsS}}
\newcommand{\BrS}{\,\mbox{\rm BrS}}
\newcommand{\TsS}{\,\mbox{\rm TsS}}

\newenvironment{lemma}[1]{\vskip+0.2cm \noindent {\bf Lemma #1.}\em }{\vskip+0.2cm}
\newenvironment{theorem}[1]{\vskip+0.2cm \noindent {\bf Theorem #1.}\em }{\vskip+0.2cm}
\newenvironment{corollary}[1]{\vskip+0.2cm \noindent {\bf Corollary #1.}\em }{\vskip+0.2cm}
\newenvironment{proposition}[1]{\vskip+0.2cm \noindent {\bf Proposition #1.}\em }{\vskip+0.2cm}

\newenvironment{definition}[1]{\vskip+0.2cm \noindent {\bf Definition #1.}}{\vskip+0.2cm}
\newenvironment{example}[1]{\vskip+0.2cm \noindent {\bf Example #1.}}{\vskip+0.2cm}

\newenvironment{remark}[1]{\vskip+0.2cm \noindent {\bf Remark #1.}}{\vskip+0.2cm}

\begin{document}

\begin{center}
    {\Large\textbf{Zero-Bidimension and Various Classes \\[0.2cm]
    of Bitopological Spaces}} \\[0.5cm]
    {\large\textbf{B. P. Dvalishvili}} \\[0.7cm]
    {\small \em Department of Algebra-Geometry of the Institute
        of Mechanics and Mathematics of the Faculty
        of Exact and Natural Sciences
        of I. Javakhishvili Tbilisi State University, \\
        2, University St., 0143 Tbilisi, Georgia}
\end{center}
\vskip+0.5cm

{\small \noindent \textbf{Abstract}
\medskip

The sum theorem and its corollaries are proved for a countable
family of zero-dimensional (in the sense of small and large
inductive bidimensions) $p\,$-closed sets, using a new notion of
relative normality whose topological correspondent is also new.
The notion of almost $n$-dimensionality is considered from the
bitopological point of view.

Bitopological spaces in which every subset is $i$-open in its
$j$-closure (i.e.,    \linebreak        $(i,j)$-submaximal spaces)
are introduced and their properties are studied. Based on the
investigations begun in [5] and [14], sufficient conditions are
found for bitopological spaces to be $(1,2)$-Baire in the class of
$p\,$-normal spaces. Furthermore, $(i,j)\dd\,I$-spaces are
introduced and both the relations between $(i,j)$-submaximal,
$(i,j)$-nodec and $(i,j)\dd\,I$-spaces, and their properties are
studied when two topologies on a set are either independent of
each other or interconnected by the inclusion, $S$-, $C$- and
$N$-relations or by their combinations.

The final part of the paper deals with the questions of
preservation of $(i,j)$-sub\-ma\-xi\-mal and $(2,1)\dd\,I$-spaces
to an image, of $D$-spaces to an image and an inverse image for
both the topological and the bitopological cases. Two theorems are
formulated containing, on the one hand, topological conditions
and, on the other hand, bitopological ones, under which a
topological space is a $D$-space.
\bigskip

\bigskip

 }

\footnotetext{\emph{E-mail address}: badri.dvalishvili@yahoo.com
(B.~Dvalishvili)}

\section*{1. Introduction}

Establishment of new results for topological as well as
bitopological spaces and strengthening of certain existing and
well-known ones have motivated this paper's systematic
investigation of different classes of bitopological spaces.

All useful notions have been collected and the following
abbreviations are used throughout the paper: $\TS$ for a
topological space, $\TsS$ for a topological subspace, $\BS$ for a
bitopological space and $\BsS$ for a bitopological subspace. The
plural form of all abbreviations is 's. Always $i,j\in \{1,2\}$,
$i\neq j$, unless stated otherwise.

Let $(X,\tau_1,\tau_2)$ be a $\BS$ and $\cP$ be some topological
property. Then $(i,j)\dd\,\cP$ denotes the analogue of this
property for $\tau_i$ with respect to $\tau_j$, and $p\,\dd\,\cP$
denotes the conjunction $(1,2)\dd\,\cP\wedge (2,1)\dd\,\cP$, that
is, $p\,\dd\,\cP$ denotes an ``absolute'' bitopological analogue
of $\cP$, where ``$p$'' is the abbreviation for ``pairwise''.
Sometimes $(1,2)\dd\,\cP\llra (2,1)\dd\,\cP$ (and thus $\llra
p\,\dd\,\cP)$ so that it suffices to consider one of these three
bitopological analogues. Moreover, there are certain cases where
equivalent topological formulations do not remain equivalent when
passing to their bitopological counterparts; in particular, this
phenomenon is observed in the case of submaximal spaces [7]. Also
note that $(X,\tau_i)$ has a property $\cP$ if and only if
$(X,\tau_1,\tau_2)$ has a property $i\dd\,\cP$, and $d\dd\,\cP$ is
equivalent to $1\dd\,\cP\wedge 2\dd\,\cP$, where ``$d$'' is the
abbreviation for ``double''.

If $\tau_1$ and $\tau_2$ are independent of each other on $X$,
then along with the properties $p\,\dd\,\cP$ and $d\dd\,\cP$ we
can cobsider the property $\sup\cP$, where $\sup\cP$ is the
$\cP$-property of the $\TS$ $(X,\sup(\tau_1,\tau_2))$, clearly not
be considered for the case $\tau_1\sbs\tau_2$, that is, in our
further notation, for a $\BS$ $(X,\tau_1<\tau_2)$.

The symbol $2^X$ is used for the power set of the set $X$, and for
a family $\cA=\{A_s\}_{s\in S}\sbs 2^X$, $\co\cA$ denotes the
conjugate family $\{X\setminus A_s:\;A_s\in \cA\}_{s\in S}$. If
$A\sbs X$, then $\tau_i\nt A$ and $\tau_i\cl A$ denote
respectively the interior and the closure of $A$ in the topology
$\tau_i$ (for a $\TS$ $(X,\tau)$ the closure of a subset $A\sbs X$
is denoted by $\ol{A})$. A set $A\sbs X$ is $p\,$-open in
$(X,\tau_1,\tau_2)$ if $A=A_1\cup A_2$, where $A_i\in \tau_i$. The
complement of a $p\,$-open set in $X$ is $p\,$-closed, that is,
$B\sbs X$ is $p\,$-closed in $(X,\tau_1,\tau_2)$ if $B=B_1\cap
B_2$, where $B_i\in\co\tau_i$ [11]. Thus, a subset $A\sbs X$ is
$p\,$-open $(p\,$-closed) in $(X,\tau_1,\tau_2)$ if and only if
$A=\tau_1\nt A\cup \tau_2\nt A$ $(A=\tau_1\cl A\cap \tau_2\cl A)$
and the family of all $p\,$-open $(p\,$-closed) subsets of a $\BS$
$(X,\tau_1,\tau_2)$ is denoted by $p\,\dd\,\cO(X)$
$(p\,\dd\Cl(X))$. It is clear that $\tau_1\cup \tau_2\sbs
p\,\dd\,\cO(X)$ $(\co\tau_1\cup\co\tau_2\sbs p\,\dd\Cl(X))$ and so
in a $\BS$ $(X,\tau_1<\tau_2)$ we have $p\,\dd\,\cO(X)=\tau_2$
$(p\,\dd\Cl(X)=\co\tau_2)$. The notion of a $p\,$-open
$(p\,$-closed) set is equivalent to the notion of a quasi open
(quasi closed) set given in [10]. The bitopological boundaries of
a subset $A\sbs X$ are $p\,$-closed sets       \linebreak
$(i,j)\dd\Fr A=\tau_i\cl A\cap \tau_j\cl (X\setminus A)$ [11].

Also, to avoid confusion with generally accepted notations, for a
$\BS$ $(X,\tau_1,\tau_2)$ we shall use the following double
indexation:
$$  A_i^d=\big\{x\in X:\;x\;\text{is an $i$-accumulation point of}
                \;A\big\}       $$
and
$$  A_j^i=\big\{x\in X:\;x\;\text{is a $j$-isolated point of}
                \;A\big\},       $$
that is, the lower indices $i$ and $j$ denote the belonging to the
topology and, therefore, $i,j\in \{1,2\}$, while the upper indices
$d$ and $i$ are fixed as the accumulation and isolation symbols,
respectively; thus $A_j^i=A\setminus A_j^d$, $A$ is a $j$-discrete
set $\llra A=A_j^i$, and
\allowdisplaybreaks
\begin{eqnarray*}
    &\ds i\dd\Bd(X)=\big\{A\in 2^X:\;\tau_i\nt A=\vnth\big\}, \\
    &\ds i\dd\,\cD(X)=\co i\dd\Bd(X)=\big\{A\in 2^X:\;\tau_i\cl A=X\big\}, \\
    &\ds i\dd\DI(X)=\{A\in 2^X:\;A\sbs A_i^d\}, \;\;\;
        (i,j)\,\dd\DI(X)=\{A\in 2^X:\;A_i^i\sbs A_j^d\}, \\
    &\ds i\dd\ST(X)=\Big\{A\in 2^X:\;A\neq\vnth,\;B\in i\dd\DI(X)\;
            \text{and}\;B\sbs A\;\text{imply}\;B=\vnth\Big\}, \\
    &\ds p\,\dd\ST(X)=\Big\{A\in 2^X:\;A\neq\vnth,\;B\in p\,\dd\DI(X)\;
            \text{and}\;B\sbs A\;\text{imply}\;B=\vnth\Big\}, \\
    &\ds i\dd\ND(X)=\big\{A\in 2^X:\;\tau_i\nt\tau_i\cl A=\vnth\big\}, \\
    &\ds (i,j)\dd\ND(X)=\big\{A\in 2^X:\;\tau_i\nt\tau_j\cl A=\vnth\big\}, \\
    &\ds (i,j)\dd\Catg_{{}_{\I}}(X)= \\
    &\ds =\Big\{A\in 2^X:\;A={\textstyle\bigcup\limits_{n=1}^\infty} A_n,\;
        A_n\in (i,j)\dd\ND(X)\;\text{for each}\;n=\ol{1,\infty}\Big\}, \\
    &\ds (i,j)\dd\Catg_{{}_{\II}}(X)=2^X\setminus (i,j)\dd\Catg_{{}_{\I}}(X), \\
    &\ds i\dd\G_\dl(X) =\Big\{A\in 2^X:\;A={\textstyle\bigcap\limits_{n=1}^\infty}
        U_n,\;U_n\in \tau_i\;\text{for each}\;n=\ol{1,\infty}\Big\}
\end{eqnarray*}
and
$$  i\dd\cF_\sg(X)=\co i\dd\G_\dl(X)=
            \Big\{A\in 2^X:\;A={\textstyle\bigcup\limits_{n=1}^\infty}
        F_n,\;F_n\in \co\tau_i\;\text{for each}\;n=\ol{1,\infty}\Big\} $$
are the families of all $i$-boundary, $i$-dense, $i$-dense in
themselves, $(i,j)$-dense in themselves, $i$-scattered,
$p\,$-scattered, $(i,j)$-first category, $(i,j)$-second category,
$i\dd\G_\dl$ and $i\dd\cF_\sg$-subsets of $X$, respectively; note
also here, that a subset $A$ of a $\BS$ $(X,\tau_1,\tau_2)$ is of
$(i,j)$-first (second) category, i.e., $A$ is of $(i,j)\dd\Catg\I$
$((i,j)\dd\Catg\II)$ if it is of $(i,j)$-first (second) category
in itself [15].

\begin{definition}{1.1}
\rm Let $(X,\tau_1,\tau_2)$ be a $\BS$. Then
\begin{enumerate}
\item[(1)] $(X,\tau_1,\tau_2)$ is $\R\dd\,p\,\dd\T_1$ (i.e.,
$p\,\dd \T_1$ in the sense of Reilly) if it is $d\dd\T_1$ [19].

\item[(2)] $(X,\tau_1,\tau_2)$ is $(i,j)$-regular if for each
point $x\in X$ and each $i$-closed set $F\sbs X$,
$x\,\ol{\in}\,F$, there exist an $i$-open set $U$ and a $j$-open
set $V$ that $x\in U$, $F\sbs V$ and $U\cap V=\vnth$ [16].

\item[(3)] $(X,\tau_1,\tau_2)$ is $p\,$-normal if for every pair
of disjoint sets $A$, $B$ in $X$, where $A$ is $1$-closed and $B$
is $2$-closed, there exist a $2$-open set $U$ and a $1$-open set
$V$ such that $A\sbs U$, $B\sbs V$and $U\cap V=\vnth$ [16].
\end{enumerate}

Moreover, $(X,\tau_1,\tau_2)$ is hereditarily $p\,$-normal if
every one of its $\BsS$ is   \linebreak      $p\,$-normal~[11].
\begin{enumerate}
\item[(4)] $(X,\tau_1,\tau_2)$ is $p\,$-connected if $X$ cannot be
expressed as a union of two disjoint sets $A$ and $B$ such that
$A\in \tau_1\setminus \{\vnth\}$ and $B\in \tau_2\setminus
\{\vnth\}$ [18] (see also [6], [8], [17]).

\item[(5)] $(X,\tau_1,\tau_2)$ is $(i,j)$-extremally disconnected
if $\tau_j\cl U=\tau_i\nt\tau_j\cl U$ for every set $U\in \tau_i$
or, equivalently, if $\tau_j\cl\tau_i\nt
A=\tau_i\nt\tau_j\cl\tau_i\nt A$ for every subset $A\sbs X$ [3].

Evidently, $(X,\tau_1,\tau_2)$ is $(1,2)$-extremally disconnected
$\llra (X,\tau_1,\tau_2)$ is       \linebreak
    $(2,1)$-extremally disconnected $\llra (X,\tau_1,\tau_2)$ is
$p\,$-extremally disconnected.

\item[(6)] $(X,\tau_1,\tau_2)$ is an $(i,j)$-Baire space or an
almost $(i,j)$-Baire space (briefly, $(i,j)\dd\BrS$ or
$A\dd(i,j)\dd\BrS)$ if every nonempty $i$-open subset of $X$ is of
\linebreak       $(i,j)$-second category or of $(i,j)$-second
category in $X$ [13], [15].

\item[(7)] $(X,\tau_1,\tau_2)$ is $(i,j)$-nodec if its every
$(i,j)$-nowhere dense subset is $j$-closed and $i$-discrete [14].
\end{enumerate}

Furthermore, in a $\BS$ $(X,\tau_1,\tau_2)$
\begin{enumerate}
\item[(8)] $\tau_1$ is coupled to $\tau_2$ (briefly,
$\tau_1C\tau_2)$ if $\tau_1\cl U\sbsq \tau_2\cl U$ for every set
$U\in \tau_1$ or, equivalently, if $\tau_1\cl\tau_1\nt A\sbsq
\tau_2\cl\tau_1\nt A$ for every subset $A\sbs X$~[22].

\item[(9)] $\tau_1$ is near $\tau_2$ (briefly, $\tau_1N\tau_2)$ if
$\tau_1\cl U\sbsq \tau_2\cl U$ for every set $U\in \tau_2$ or,
equivalently, if $\tau_1\cl\tau_2\nt A\sbsq \tau_2\cl\tau_2\nt A$
for every subset $A\sbs X$ [15].

\item[(10)] $\tau_1$ and $\tau_2$ are $S$-related on $X$ (briefly,
$\tau_1S\tau_2)$ if $\tau_1\nt A\sbs \tau_1\cl\tau_2\nt A\wedge
\tau_2\nt A\sbs \tau_2\cl\tau_1\nt A$ for every subset $A\sbs X$
[20].
\end{enumerate}

Since for a $\BS$ $(X,\tau_1<\tau_2)$ the inclusions
$$  \begin{matrix}
        (2,1)\dd\ND(X) &\!\!\sbs\!\!& 2\dd\ND(X) \\
            \cap  & & \cap \\
        1\dd\ND(X) &\!\!\sbs\!\!& (1,2)\dd\ND(X)
    \end{matrix}        $$
are correct, in the case where $\tau_1\sbs \tau_2$ we come to the
following evident implications:
$$  (X,\tau_1,\tau_2)\;\text{is $1$-nodec} \lra
        (X,\tau_1,\tau_2)\;\text{is $(2,1)$-nodec}      $$
and
$$  (X,\tau_1,\tau_2)\;\text{is $(1,2)$-nodec} \lra
        (X,\tau_1,\tau_2)\;\text{is $2$-nodec}\,.     $$

Moreover, according to (1) of Theorem~2.1.10 in [15], for a $\BS$
$(X,\tau_1<_S\tau_2)$, where $\tau_1<_S\tau_2\llra
(\tau_1\sbs\tau_2\wedge\tau_1S\tau_2)$, in addition to the above
implications, we have:
$$  \begin{matrix}
        (X,\tau_1,\tau_2)\;\text{is $1$-nodec} &\!\!\lra\!\!&
            (X,\tau_1,\tau_2)\;\text{is $(2,1)$-nodec} \\
        \Downarrow && \Downarrow \\
        (X,\tau_1,\tau_2)\;\text{is $(1,2)$-nodec} &\!\!\lra\!\!&
            (X,\tau_1,\tau_2)\;\text{is $2$-nodec}\,.
    \end{matrix}        $$
\end{definition}

\begin{remark}{1.2}
\rm Every $\BsS$ $(Y,\tau_1',\tau_2')$ of a $(1,2)$-nodec $\BS$
$(X,\tau_1<\tau_2)$ is also $(1,2)$-nodec. Indeed: if $A\in
(1,2)\dd\ND(Y)$, then by (1) of Theorem~1.5.13 in [15], $A\in
(1,2)\dd\ND(X)$ and, hence, $A=\tau_2\cl A=A_1^i$. Thus $A$ is
$2$-closed and $1$-discrete in $(Y,\tau_1',\tau_2')$ too.
\end{remark}

\begin{definition}{1.3}
\rm A subset $A$ of a $\BS$ $(X,\tau_1,\tau_2)$ is said to be
$(i,j)$-locally closed if for each point $x\in A$ there exists a
set $U\in \tau_i$ such that $x\in U$ and $U\cap A=U\cap\tau_j\cl
A$ [15].
\end{definition}

The families of all such subsets of $X$ are denoted by
$(i,j)\dd\LC(X)$ and it is not difficult to see that
\begin{eqnarray*}
    &\ds A\in (i,j)\dd\LC(X)\llra \big(A\in \tau_i' \;\;\mbox{in the}\;\;
            \BsS \;\; (\tau_j\cl A,\tau_1',\tau_2')\big)\llra \\
    &\ds \llra \big(A=U\cap F, \;\;\mbox{where}\;\;
            U\in \tau_i \;\;\mbox{and}\;\;F\in \co\tau_j\big).
\end{eqnarray*}

Hence $\tau_i\cup\co\tau_j\sbs (i,j)\dd\LC(X)$ and for a $\BS$
$(X,\tau_1<\tau_2)$ the following inclusions hold:
$$  \begin{matrix}
        \tau_2\cup\co\tau_2 &\!\!\sbs\!\!& 2\dd\LC(X)
            &\!\!\sps\!\!& (2,1)\dd\LC(X) &\!\!\sps\!\!&
                \tau_2\cup\co\tau_1 \\
        & & \cup & & \cup & & \\
        \tau_1\cup\co\tau_2 &\!\!\sbs\!\!& (1,2)\dd\LC(X)
            &\!\!\sps\!\!& 1\dd\LC(X) &\!\!\sps\!\!&
                \tau_1\cup\co\tau_1\,.
    \end{matrix}        $$

\begin{definition}{1.4}
\rm Let $(x,A)$ be a pair in a $\BS$ $(X,\tau_1,\tau_2)$ such that
$A\in \co\tau_i$ and $x\,\ol{\in}\,A$. Then a $p\,$-closed set $T$
is a partition, corresponding to the pair $(x,A)$, if $X\setminus
T=H_1\cup H_2$, where $H_i\in \tau_i\setminus \{\vnth\}$, $x\in
H_i$, $A\sbs H_j$ and $H_1\cap H_2=\vnth$.
\end{definition}

If $(x,A)$ is a pair in the above sense, then one can easily
verify that for an $i$-open neighborhood $U(x)$ $(j$-open
neighborhood $U(A))$ such that $\tau_j\cl U(x)\sbs X\setminus A$
$(\tau_i\cl U(A)\sbs X\setminus \{x\})$, the sets $(j,i)\dd\Fr
U(x)$ $((i,j)\dd\Fr U(A))$ are partitions, corresponding to
$(x,A)$ and, conversely, if $T$ is a partition, corresponding to
$(x,A)$, then $(j,i)\dd\Fr H_i\sbs T$.

For the pairwise small inductive dimension $p\,\dd\ind X$ we have:

$p\,\dd\ind X=0\llra (\tau_1$ has a base consisting of $2$-closed
sets and $\tau_2$ has a base consisting of $1$-closed sets)
$\llra$ (the empty set is a partition corresponding to any pair
$(x,A)$, where $A\in \co\tau_1$, $x\,\ol{\in}\,A$ and the empty
set is a partition corresponding to any pair $(x,A)$, where $A\in
\co\tau_2$, $x\,\ol{\in}\,A)$ [11], [12].

\begin{definition}{1.5}
\rm Let $(A,B)$ be a pair of subsets of a $\BS$
$(X,\tau_1,\tau_2)$ such that $A\in \co\tau_1$, $B\in \co\tau_2$
and $A\cap B=\vnth$. Then we say that a $p\,$-closed set $T$ is a
partition, corresponding to $(A,B)$, if $X\setminus T=H_1\cup
H_2$, where $H_i\in \tau_i\setminus\{\vnth\}$, $A\sbs H_2$, $B\sbs
H_1$ and $H_1\cap H_2=\vnth$.
\end{definition}

If $(A,B)$ is a pair in the above sense and there exists a
$2$-open neighborhood $U(A)$ $(1$-open neighborhood $U(B))$ such
that $\tau_1\cl U(A)\sbs X\setminus B$ $(\tau_2\cl U(B)\sbs
X\setminus A)$, then $(1,2)\dd\Fr U(A)$ $((2,1)\dd\Fr U(B))$ is a
partition, corresponding to $(A,B)$ and, conversely, if $T$ is a
partition, corresponding to $(A,B)$, then $(j,i)\dd\Fr H_i\sbs T$.

Now, for the pairwise large inductive dimension $p\,\dd\Ind X$ we
have:

$p\,\dd\Ind X=0\llra$ (the empty set is a partition corresponding
to any pair $(A,B)$, where $A\in \co\tau_1$, $B\in \co\tau_2$ and
$A\cap B=\vnth\;)\llra$ (for every $1$-closed set $F$ and any
$2$-neighborhood $U(F)$ there exists a neighborhood
$V(F)\in\tau_2\cap\co\tau_1$ such that $V(F)\sbs U(F)$ and, for
every $2$-closed set $\Phi$ and any $1$-neighborhood $U(\Phi)$
there exists a neighborhood $V(\Phi)\in\tau_1\cap\co\tau_2$ such
that $V(\Phi)\sbs U(\Phi))$ [11], [12].

\begin{definition}{1.6}
\rm Let $(X,\tau_1,\tau_2)$ and $(Y,\gm_1,\gm_2)$ be $\BS$'s. Then
a function   \linebreak       $f: (X,\tau_1,\tau_2)\to
(Y,\gm_1,\gm_2)$ is said to be $i$-continuous $(i$-open,
$i$-closed) if the induced functions $f: (X,\tau_i)\to (Y,\gm_i)$
are continuous (open, closed) [18].
\end{definition}

Finally, please note that all bitopological generalizations are
constructed in the commonly accepted manner, so that if the
topologies coincide, we obtain the original topological notions.

\section*{2. Some Special Notions and the Sum Theorem for Bidimension Zero}

Our first step in this section is to prove some simple, but
nevertheless principal facts, having an independent interest.

\begin{proposition}{2.1}
For any subsets $A\in\tau_2$ $(A\in \co\tau_1)$ and $B\in\tau_1$
$(B\in \co\tau_2)$ of a $\BS$ $(X,\tau_1,\tau_2)$ the sets
$A\setminus B$ and $B\setminus A$ are $p\,$-separated, that is,
$$  \big(\tau_1\cl(A\setminus B)\cap (B\setminus A)\big)\cup
        \big((A\setminus B)\cap \tau_2\cl(B\setminus A)\big)=\vnth.  $$
\end{proposition}

\noindent \textbf{Proof.\/} We will prove the condition, given
without brackets, since the other one can be proved similarly.

Evidently, $A=A\cap (A\cup B)\in \tau_2'$ and $B=B\cap (A\cup
B)\in \tau_1'$ in $(A\cup B,\tau_1',\tau_2')$. Therefore
$A\setminus B=(A\cup B)\setminus B\in \co\tau_1'$ and $B\setminus
A=(A\cup B)\setminus A\in \co\tau_2'$. Let us consider the $\BsS$
$((A\setminus B)\cup (B\setminus A),\tau_1'',\tau_2'')$ of the
$\BS$ $(A\cup B,\tau_1',\tau_2')$. Then
$$  (A\setminus B)\cap \big((A\setminus B)\cup (B\setminus A)\big)=
        A\setminus B\in\co\tau_1''      $$
and
$$  (B\setminus A)\cap \big((A\setminus B)\cup (B\setminus A)\big)=
        B\setminus A\in\co\tau_2''.      $$
Hence $A\setminus B\in\tau_2''$ and $B\setminus A\in\tau_1''$ so
that
$$  \big(\tau_1''\cl(A\setminus B)\cap (B\setminus A)\big)\cup
        \big((A\setminus B)\cap \tau_2''\cl(B\setminus A)\big)=\vnth $$
and, thus,
$$   \big(\tau_1\cl(A\setminus B)\cap (B\setminus A)\big)\cup
        \big((A\setminus B)\cap \tau_2\cl(B\setminus A)\big)=\vnth.
                    \;\; \Box   $$

\begin{definition}{2.2}
\rm Let $(Y,\tau_1',\tau_2')$ and $(Z,\tau_1'',\tau_2'')$ be
$\BsS$'s of a $\BS$ $(X,\tau_1,\tau_2)$ such that $X=Y\cup Z$.
Then the $(i,j)$-closures of any subset $A\sbs X$, corresponding
to the pair $(Y,Z)$, are the sets
$$  (i,j)\dd\cl A(Y,Z)=
        \tau_i\cl(A\cap Y)\cup\tau_j\cl\big(A\cap (Z\setminus Y)\big). $$
Moreover,
$$  D\dd\cl A(Y,Z)=\big((1,2)\dd\cl A(Y,Z)\cap Y\big)\cup
        \big((2,1)\dd\cl A(Z,Y)\cap Z\big)      $$
and $A$ is $D$-closed with respect to $(Y,Z)$ if $A=D\dd\cl
A(Y,Z)$.
\end{definition}

\begin{proposition}{2.3}
If for a $\BS$ $(X,\tau_1,\tau_2)$ there are $\BsS$'s
$(Y,\tau_1',\tau_2')$, $(Z,\tau_1'',\tau_2'')$ such that $X=Y\cup
Z$ and the sets $Y\setminus Z$, $Z\setminus Y$ are
$p\,$-separated, then for every subset $A\sbs X$ we have
$$  D\dd\cl A(Y,Z)=\tau_1'\cl(A\cap Y)\cup\tau_2''\cl(A\cap Z). $$
\end{proposition}

\noindent \textbf{Proof.\/} Since $\tau_2\cl(Z\setminus Y)\cap
(Y\setminus Z)=\vnth$, we have $\tau_2\cl(Z\setminus Y)\sbs Z$ so
that
$$  \tau_2\cl\big(A\cap (Z\setminus Y)\big)\cap Y\sbs
        \tau_2\cl\big(A\cap (Z\setminus Y)\big)\sbs
            \tau_2\cl(A\cap Z)\cap Z.       $$
Hence
$$  \big(\tau_1\cl(A\cap Y)\cup
            \tau_2\cl(A\cap (Z\setminus Y))\big)\cap Y=
        (\tau_1\cl(A\cap Y)\cap Y\big)\cup M,       $$
where
$$  M=\tau_2\cl\big(A\cap (Z\setminus Y)\big)\cap Y\sbs
            \tau_2\cl (A\cap Z)\cap Z.        $$

By the similar manner we obtain that
$$  \big(\tau_2\cl(A\cap Z)\cup
            \tau_1\cl(A\cap (Y\setminus Z))\big)\cap Z=
        \big(\tau_2\cl (A\cap Z)\cap Z\big)\cup N,      $$
where
$$  N=\tau_1\cl (A\cap (Y\setminus Z))\cap Z\sbs
        \tau_1\cl (A\cap Y)\cap Y.      $$
Thus
\begin{eqnarray*}
    &\ds \Big(\big(\tau_1\cl (A\cap Y)\cup
        \tau_2\cl(A\cap (Z\setminus Y))\big)\cap Y\Big)\cup \\
    &\ds \cup \Big(\big(\tau_2\cl (A\cap Z)\cup
        \tau_1\cl(A\cap (Y\setminus Z))\big)\cap Z\Big)= \\
    &\ds =\big(\tau_1\cl (A\cap Y)\cap Y\big)\cup
            \big(\tau_2\cl(A\cap Z)\cap Z\big)=
        \tau_1'\cl(A\cap Y)\cap \tau_2''\cl(A\cap Z). \;\; \Box
\end{eqnarray*}

\begin{corollary}{2.4}
If $(Y,\tau_1',\tau_2')$ and $(Z,\tau_1'',\tau_2'')$ are $\BsS$'s
of a $\BS$ $(X,\tau_1,\tau_2)$, where $Y\in \tau_2$ $(Y\in
\co\tau_1)$, $Z\in \tau_1$ $(Z\in \co\tau_2)$ and $X=Y\cup Z$,
then
$$  D\dd\cl A(Y,Z)=\tau_1'\cl(A\cap Y)\cup \tau_2''\cl(A\cap Z). $$
\end{corollary}

\noindent \textbf{Proof.\/} By Proposition~2.1, the sets
$Y\setminus Z$ and $Z\setminus Y$ are $p\,$-separated and, hence,
it remains to use Proposition~2.3.~$\Box$

\begin{corollary}{2.5}
If $(Y,\tau_1',\tau_2')$ and $(Z,\tau_1'',\tau_2'')$ are $\BsS$'s
of a $\BS$ $(X,\tau_1,\tau_2)$, $X=Y\cup Z$ and $Y\setminus Z$,
$Z\setminus Y$ are $p\,$-separated, then $A=D\dd\cl A(Y,Z)$ if and
only if $A\cap Y\in\co\tau_1'$ and $A\cap Z\in \co\tau_2''$.
\end{corollary}

\noindent \textbf{Proof.\/} If $A=D\dd\cl A(Y,Z)$, then according
to Proposition~2.3,
$$  A=\tau_1'\cl(A\cap Y)\cup \tau_2''\cl(A\cap Z)=
        (A\cap Y)\cup (A\cap Z)   $$
and so $A\cap Y\in \co\tau_1'$, $A\cap Z\in \co\tau_2''$.
Conversely, if $A\cap Y\in \co\tau_1'$ and $A\cap
Z\in\co\tau_2''$, then once more applying Proposition~2.3, we
obtain that
\begin{eqnarray*}
    &\ds D\dd\cl A(Y,Z)=\tau_1'\cl(A\cap Y)\cup \tau_2''\cl(A\cap Z)= \\
    &\ds =(A\cap Y)\cup (A\cap Z)=A\cap (Y\cup Z)=A\cap X=A. \;\; \Box
\end{eqnarray*}

\begin{corollary}{2.6}
If $(Y,\tau_1',\tau_2')$ and $(Z,\tau_1'',\tau_2'')$ are $\BsS$'s
of a $\BS$ $(X,\tau_1,\tau_2)$ such that $X=Y\cup Z$. Then for
every subset $A\sbs X$ the following conditions hold:
\begin{enumerate}
\item[$(1)$] $A\sbs (1,2)\dd\cl A(Y,Z)\cap (2,1)\dd\cl A(Z,Y)$.

\item[$(2)$] $(1,2)\dd\cl A(Y,Z)\cup (2,1)\dd\cl
A(Z,Y)=\tau_1\cl(A\cap Y)\cup\tau_2\cl(A\cap Z)$.

\item[$(3)$] If $Y\cap Z=\vnth$, then
$$  (1,2)\dd\cl A(Y,Z)=(2,1)\dd\cl A(Z,Y)=
        \tau_1\cl(A\cap Y)\cup \tau_2\cl(A\cap Z).  $$

\item[$(4)$] If $Y=\vnth$, that is, if $X=Z$, then
$$  (1,2)\dd\cl A(\vnth,Z)=\tau_2\cl A=(2,1)\dd\cl A(Z,\vnth). $$

\item[$(5)$] If $Z=\vnth$, that is, if $X=Y$, then
$$  (1,2)\dd\cl A(Y,\vnth)=\tau_1\cl A=(2,1)\dd\cl A(\vnth,Y). $$

\item[$(6)$] If $A\sbs Y$, then $(1,2)\dd\cl A(Y,Z)=\tau_1\cl A$
and if $A\sbs Y\setminus Z$, then
$$  (1,2)\dd\cl A(Y,Z)=\tau_1\cl A=(2,1)\dd\cl A(Z,Y). $$

\item[$(7)$] If $A\sbs Z$, then $(2,1)\dd\cl A(Z,Y)=\tau_2\cl A$
and if $A\sbs Z\setminus Y$, then
$$  (1,2)\dd\cl A(Y,Z)=\tau_2\cl A=(2,1)\dd\cl A(Z,Y). $$

\item[$(8)$] If $Y=Z=X$, then
$$  (1,2)\dd\cl A(X,X)=\tau_1\cl A, \;\;\;
        (2,1)\dd\cl A(X,X)=\tau_2\cl A      $$
and
$$  D\dd\cl A(X,X)=\tau_1\cl A\cup \tau_2\cl A.   $$

\item[$(9)$] If $Y=A$, $Z=X\setminus A$, then
$$  (i,j)\dd\Fr A=(i,j)\dd\cl A(A,X\setminus A)\cap
            (j,i)\dd\cl(X\setminus A)(X\setminus A,A).      $$
\end{enumerate}

Moreover, if $\tau_1\sbs \tau_2$, then
\begin{enumerate}
\item[$(10)$] $\tau_2\cl A\sbs (1,2)\dd\cl A(Y,Z)\cap (2,1)\dd\cl
A(Z,Y)\sbs \tau_1\cl A$.

\item[$(11)$] If $A\sbs Y$, then $\tau_2\cl A\sbs (2,1)\dd\cl
A(Z,Y)\sbs (1,2)\cl A(Y,Z)=\tau_1\cl A$.

\item[$(12)$] If $A\sbs Z$, then $\tau_2\cl A=(2,1)\dd\cl
A(Z,Y)\sbs (1,2)\dd\cl A(Y,Z)\sbs \tau_1\cl A$.

\item[$(13)$] $D\dd\cl A(X,X)=\tau_1\cl A$.
\end{enumerate}
\end{corollary}

The proof of (1)--(13) is trivial.

Clearly, the (i,j)-closures and $D$-closure have many other
interesting properties.

\begin{proposition}{2.7}
For a $\BS$ $(X,\tau_1,\tau_2)$ the following conditions are
equivalent:
\begin{enumerate}
\item[$(1)$] $(X,\tau_1,\tau_2)$ is hereditarily $p\,$-normal.

\item[$(2)$] Every $p\,$-open subset of $X$ is $p\,$-normal.

\item[$(3)$] If $A$, $B$ are $p\,$-separated subsets of $X$, then
there are sets $U\in \tau_2$, $V\in\tau_1$ such that $A\sbs U$,
$B\sbs V$ and $U\cap V=\vnth$.
\end{enumerate}
\end{proposition}

\noindent \textbf{Proof.\/} The equivalence $(1)\llra (3)$ is
proved in Theorem~0.2.2 from [15]. The implication $(1)\lra (2)$
is obvious. The proof of the implication $(2)\lra (3)$ is given in
the first part of the proof of Theorem~0.2.2 from [15], since the
set $Y=X\setminus (\tau_1\cl A\cap \tau_2\cl B)$ is $p\,$-open.
$\Box$

\begin{corollary}{2.8}
If a $\BS$ $(X,\tau<_C\tau_2)$ is hereditarily $p\,$-normal, then
it is hereditarily $1$-normal.
\end{corollary}

\noindent \textbf{Proof.\/} By the well-known topological fact it
is sufficient to prove only that every $1$-open $\BsS$ of $X$ is
$1$-normal. Let $U\in \tau_1\setminus \{\vnth\}$ be any set. Then
$U\in p\,\dd\,\cO(X)$ and by (2) of Proposition~2.7,
$(U,\tau_1',\tau_2')$ is $p\,$-normal. On the other hand, by (2)
of Corollary~2.2.8 in [15], $\tau_1'<_C\tau_2'$ and, hence, by (4)
of the same corollary, $(U,\tau_1',\tau_2')$ is $1$-normal.~$\Box$

\begin{corollary}{2.9}
If a $\BS$ $(X,\tau<_N\tau_2)$ is hereditarily $2$-normal, then it
is hereditarily $p\,$-normal and so it is hereditarily $1$-normal.
\end{corollary}

\noindent \textbf{Proof.\/} When $\tau_1<_N\tau_2$, then also
$\tau_1<_C\tau_2$ (Corollary~2.3.10 in [15]) and by Corollary~2.8
it suffices to prove only that every hereditarily $2$-normal $\BS$
\linebreak     $(X,\tau<_N\tau_2)$ is hereditarily $p\,$-normal.
Let $U\in p\,\dd\,\cO(X)$ be any set. Then $U\in \tau_2$ as
$\tau_1\sbs\tau_2$, so that $(U,\tau_1',\tau_2')$ is $2$-normal.
But, according to (2) of Corollary~2.2.13 in [15],
$\tau_1'<_N\tau_2'$ and, hence, by (4) of the same corollary,
$(U,\tau_1',\tau_2')$ is $p\,$-normal.~$\Box$

\begin{definition}{2.10}
\rm A $\BS$ $(X,\tau_1,\tau_2)$ is $p\,$-hypernormal if for every
pair $A$, $B$ of $p\,$-separated subsets of $X$ (i.e., $(\tau_1\cl
A\cap B)\cup (A\cap \tau_2\cl B)=\vnth)$ there are sets $U\in
\tau_2$, $V\in \tau_1$ such that $A\sbs U$, $B\sbs V$ and
$\tau_1\cl U\cap\tau_2\cl V=\vnth$.
\end{definition}

Take place the following

\begin{proposition}{2.11}
A $\BS$ $(X,\tau_1,\tau_2)$ is $p\,$-hypernormal if and only if it
is        \linebreak       $p\,$-extremally disconnected and
hereditarily $p\,$-normal.
\end{proposition}

\noindent \textbf{Proof.\/} Evidently, every $p\,$-extremally
disconnected and hereditarily $p\,$-nor\-mal $\BS$ is
$p\,$-hypernormal. Therefore, suppose that $(X,\tau_1,\tau_2)$ is
$p\,$-hy\-per\-nor\-mal. Then, by Proposition~2.7,
$(X,\tau_1,\tau_2)$ is hereditarily $p\,$-normal and so, by (5) of
Definition~1.1, it remains to prove only that $\tau_j\cl U\in
\tau_i$ for every $i$-open subset $U\sbs X$. Since $U\in \tau_i$,
we have
$$  \big(U\cap \tau_i\cl(X\setminus \tau_j\cl U)\big)\cup
        \big(\tau_j\cl U\cap (X\setminus \tau_j\cl U)\big)=\vnth $$
so that $U$ and $X\setminus \tau_j\cl U$ are $p\,$-separated and,
hence, by condition, there are sets $M\in \tau_i$, $N\in \tau_j$
such that $U\sbs M$, $(X\setminus \tau_j\cl U)\sbs N$ and
$\tau_j\cl M\cap \tau_i\cl N=\vnth$. Since $\tau_j\cl U\sbs
\tau_j\cl M$ and $\tau_i\cl(X\setminus \tau_j\cl U)\sbs \tau_i\cl
N$, we obtain that $\tau_j\cl U\cap \tau_i\cl(X\setminus \tau_j\cl
U)=\vnth$ so that
$$  \tau_j\cl U\cap (X\setminus \tau_i\nt\tau_j\cl U)=\vnth $$
and, thus, $\tau_j\cl U=\tau_i\nt\tau_j\cl U$.~$\Box$

\begin{corollary}{2.12}
If $(X,\tau_1,\tau_2)$ is $p\,$-hypernormal, then
\begin{enumerate}
\item[$(1)$] $p\,\dd\Ind X=0$ or, equivalently, $p\,\dd\,\dim
X=0$.

\item[$(2)$] If $(X,\tau_1,\tau_2)$ is $\R\dd\,p\,\dd\T_1$, then
$$  p\,\dd\ind X=p\,\dd\Ind X=p\,\dd\,\dim X=0.  $$
\end{enumerate}

Moreover, in addition to the above conditions, for a $\BS$
$(X,\tau_1<_C\tau_2)$ we have:
\begin{enumerate}
\item[$(3)$] $(X,\tau_1,\tau_2)$ is $1$-hypernormal.

\item[$(4)$] $1\dd\Ind X=1\dd\,\dim X=p\,\dd\,\Ind X=p\,\dd\,\dim
X=0$ \\
and if $(X,\tau_1,\tau_2)$ is $1\dd\T_1$, then
$$  1\dd\Ind X=1\dd\,\dim X=p\,\dd\Ind X=p\,\dd\,\dim X=d\,\dd\ind X=
        p\,\dd\ind X=0.    $$
\end{enumerate}
\end{corollary}

\noindent \textbf{Proof.\/} (1) The equality $p\,\dd\Ind X=0$ is
obvious, since, by Proposition~2.11, $(X,\tau_1,\tau_2)$ is
$p\,$-normal and $p\,$-extremally disconnected. Hence, it remains
to use (1) of Corollary~3.3.5 in [15].

(2) It suffices to use (2) of Corollary 3.3.5 in [15].

(3) According to Corollary 2.8, $(X,\tau_1,\tau_2)$ is
hereditarily $1$-normal, and by (8) of Corollary~2.2.8 in [15],
$(X,\tau_1,\tau_2)$ is $1$-extremally disconnected.

(4) The first part follows directly from the well-known
topological fact, (1) above and (1) of Theorem~3.2.38 in [15]. The
rest is an immediate consequence of (2) above and (5) of
Theorem~3.1.36 in [15].

\begin{corollary}{2.13}
For a $\BS$ $(X,\tau_1<_N\tau_2)$ the following implications hold:
\begin{eqnarray*}
    &&\hskip-0.6cm (X,\tau_1,\tau_2)\;\mbox{is $d$-hypernormal} \lla
        (X,\tau_1,\tau_2)\;\mbox{is $2$-hypernormal} \lra \\
    \lra &&\hskip-0.6cm (X,\tau_1,\tau_2)\;\mbox{is $p\,$-hypernormal} \lra
        (X,\tau_1,\tau_2)\;\mbox{is $1$-hypernormal}.
\end{eqnarray*}
\end{corollary}

\noindent \textbf{Proof.\/} If $(X,\tau_1,\tau_2)$ is
$2$-hypernormal, then by Corollary~2.9, it is hereditarily
$p\,$-normal, and by (7) of Corollary~2.3.13 in [15], it is also
$p\,$-extremally disconnected. Hence, according to
Proposition~2.11, $(X,\tau_1,\tau_2)$ is $p\,$-hypernormal.
Furthermore, by Corollary~2.3.10 in [15], $\tau_1<_C\tau_2$ and
so, by (3) of Corollary~2.12, $(X,\tau_1,\tau_2)$ is
$1$-hypernormal.

Now, the first implication from right to left is obvious.~$\Box$
\vskip+0.2cm

The rest of this section is devoted to new notions of relative
normality of $\BsS$'s and their applications in the theory of
dimension of $\BS$'s. Their topological counterpart is new as
well. As will be shown below, these notions prove to be the key
tool for correcting an error made in proving Theorem~3.2.26 and
its Co\-rol\-la\-ri\-es~3.2.27--3.2.30 in [15]. In particular, it
allows us to prove the above-mentioned theorem and its corollaries
for a sequence of $p\,$-closed sets. The latter circumstance
emphasizes once more a special role which relative (bi)topological
properties play not only in the development of respective
theories, but also in the strengthening of the previously known
results.

\begin{definition}{2.14}
\rm We will say that a $\TsS$ $(Y,\tau')$ of a $\TS$ $(X,\tau)$ is
$\WS$-supernormal in $X$ if for each pair of disjoint sets $A$,
$B$, where $A$ is closed in $(Y,\tau')$ and $B$ is closed in
$(X,\tau)$, there are disjoint sets $U$, $V$ such that $U$ is open
in $(Y,\tau')$, $V$ is open in $(X,\tau)$ and $A\sbs U$, $B\sbs
V$.
\end{definition}

It is clear that if $(X,\tau)$ is normal, then every closed
subspace of $(X,\tau)$ is $WS$-normal in $X$, and if $(Y,\tau')$
is $\WS$-supernormal in $X$, then $(Y,\tau')$ is normal (in
itself). But for the opposite implication we have the following
elementary

\begin{example}{2.15}
\rm Let us consider the natural $\TS$ $(\bR,\om)$ and an open
interval $(a,b)\sbs \bR$, which is normal (in itself). If
$A=[c,b)$ and $B=[b,d]$ are sets, where $a<c<b<d$, then  $A\in
\co\om'$ in $((a,b),\om')$, $B\in \co\om$ and $A\cap B=\vnth$. But
for any neighborhood $V(B)\in \om$ we have $A\cap V(B)\neq\vnth$.
\end{example}

Now, we give the bitopological modifications of relative
$\WS$-supernormality in the above sense and of relative strong
normality in the sense of [4].

\pagebreak
\begin{definition}{2.16}
\rm We will say that a $\BsS$ $(Y,\tau_1',\tau_2')$ of a $\BS$
$(X,\tau_1,\tau_2)$ is $(i,j)\dd\WS$-supernormal $(p\,$-strongly
normal) in $X$ if for each pair of disjoint sets $A\in
\co\tau_i'$, $B\in \co\tau_j$ $(A\in \co\tau_1'$, $B\in
\co\tau_2')$ there are disjoint sets $U\in \tau_j'$, $V\in \tau_i$
$(U\in \tau_2$, $V\in \tau_1)$ such that $A\sbs U$ and $B\sbs V$.
\end{definition}

Evidently, if $(X,\tau_1,\tau_2)$ is $p\,$-normal, then for every
subset $Y\in \co\tau_i$ the $\BsS$ $(Y,\tau_1',\tau_2')$ is
$(i,j)\dd\WS$-supernormal in $X$ and so, if $(X,\tau_1,\tau_2)$ is
$p\,$-normal and $Y\in \co\tau_1\cap \co\tau_2$, then
$(Y,\tau_1',\tau_2')$ is $p\,\dd\WS$-supernormal in $X$ (clearly,
for a $\BS$ $(X,\tau_1<\tau_2)$ the last fact is correct for $Y\in
\co\tau_1)$.

\begin{example}{2.17}
\rm Let $(\bR,\om_1,\om_2)$ be the natural $\BS$, that is,
$\om_1=\{\vnth,\bR\}\cup\{(a,+\infty):\;a\in \bR\}$ and
$\om_2=\{\vnth,\bR\}\cup\{(-\infty,a):\;a\in \bR\}$ are the lower
and upper topologies, respectively. Then it is not difficult to
see that every set $Y\in \om_i\setminus \{\vnth\}$ is
$(i,j)\dd\WS$-supernormal in the $p\,$-normal $\BS$
$(\bR,\om_1,\om_2)$, but it is not $(j,i)\dd\WS$-supernormal in
$(\bR,\om_1,\om_2)$. Moreover, if $Y=[a,b]=(-\infty,b]\cap
[a,+\infty)\in p\,\dd\Cl(\bR)$, then it is $p\,\dd\WS$-supernormal
in $(\bR,\om_1,\om_2)$.
\end{example}

\begin{proposition}{2.18}
In a $\BS$ $(X,\tau_1<\tau_2)$ we have:
\begin{enumerate}
\item[$(1)$] If $Y\in \tau_2$ and $\tau_1<_N\tau_2$, then the
following implications hold:
$$ \ \hskip-1cm \begin{matrix}
        \text{$(Y,\tau_1',\tau_2')$\,is\,$2\dd\WS$-supernorm.\,in\,$X$}
                    \!\!\!&\!\!\!\lra\!\!\!&\!\!\!
            \text{$(Y,\tau_1',\tau_2')$\,is\,$(1,2)\dd\WS$-supernorm.\,in\,$X$} \\
        \big\Downarrow & & \big\Downarrow \\
        \text{$(Y,\tau_1',\tau_2')$\,is\,$(2,1)\dd\WS$-supernorm.\,in\,$X$}
                    \!\!\!&\!\!\!\lra\!\!\!&\!\!\!
            \text{$(Y,\tau_1',\tau_2')$\,is\,$1\dd\WS$-supernorm.\,in\,$X$}\,.
    \end{matrix}        $$

\item[$(2)$] If $\tau_1<_C\tau_2$, then
\begin{eqnarray*}
    &\ds \ \hskip-0.7cm (Y,\tau_1',\tau_2') \;\mbox{is $p\,$-strongly normal in}\;
            X \lra (Y,\tau_1',\tau_2') \;\mbox{is $1$-strongly
                    normal in}\; X
\end{eqnarray*}
and if $\tau_1<_N\tau_2$, then
$$  \xymatrix{
        Y\;\mbox{is $2$-strongly normal in} \ar@{=>}[r] \ar@{=>}[d]
            & Y\;\mbox{is $d$-strongly normal in}\;X \ar@{=>}[d] \\
        Y\;\mbox{is $p\,$-strongly normal in}\ar@{=>}[r]
            & Y\;\mbox{is $1$-strongly normal in} }     $$
\end{enumerate}
\end{proposition}

\noindent \textbf{Proof.\/} (1) First of all, let us note that if
$Y\in \tau_2$, then by (2) of Corollary~2.3.13 in [15],
$\tau_1<_N\tau_2$ implies $\tau_1'<_N\tau_2'$.

For the upper horizontal implication let $A\in \co\tau_1'$, $B\in
\co\tau_2$ and $A\cap B=\vnth$. Since $A\in \co\tau_1'\sbs
\co\tau_2'$ and $(Y,\tau_1',\tau_2')$ is $2\dd\WS$-supernormal in
$X$, there are neighborhoods $U'(A)\in \tau_2'\sbs \tau_2$,
$U'(B)\in \tau_2$ such that $U'(A)\cap U'(B)=\vnth$ and so
$\tau_2\cl U'(A)\cap U'(B)=\vnth$. According to (2) of
Corollary~2.3.12 in [15],
$$  U'(A)\sbs \tau_2\nt\tau_2\cl U'(A)=\tau_1\nt\tau_2\cl U'(A)=
        U''(A)\in \tau_1\sbs\tau_2  $$
and for the set $U(A)=U''(A)\cap Y\in \tau_1'\sbs \tau_2'$ we have
$U(A)\cap \tau_1\cl U'(B)=\vnth$. Moreover, by Corollary~2.3.10 in
[15],
$$  \tau_1<_N\tau_2\lra \tau_1<_C\tau_2      $$
and (2) of Corollary~2.2.7 in [15] gives that
$$  U'(B)\sbs \tau_2\nt\tau_1\cl U'(B)=\tau_1\nt\tau_1\cl U'(B)=
        U(B)\in \tau_1.     $$
Now, it is clear that $U(A)\cap U(B)=\vnth$ and, hence,
$(Y,\tau_1',\tau_2')$ is $(2,1)\dd\WS$-supernormal in $X$.

For the lower horizontal implication let $A\in \co\tau_1'$, $B\in
\co\tau_1$ and $A\cap B=\vnth$. Since $A\in \co\tau_1'\sbs
\co\tau_2'$ and $(Y,\tau_1',\tau_2')$ is $(2,1)\dd\WS$-supernormal
in $X$, there are neighborhoods $U(A)\in \tau_1'$, $U'(B)\in
\tau_2$ such that $U(A)\cap U'(B)=\vnth$. Since $U(A)\in
\tau_1'\sbs\tau_2'\sbs \tau_2$, we have $U(A)\cap \tau_2\cl
U'(B)=\vnth$. But
$$  U'(B)\sbs \tau_2\nt\tau_2\cl U'(B)=\tau_1\nt\tau_2\cl U'(B)=
            U(B)\in \tau_1      $$
and $U(A)\cap U(B)=\vnth$. Thus $(Y,\tau_1',\tau_2')$ is
$1\dd\WS$-supernormal in $X$.

For the left vertical implication let $A\in \co\tau_2'$, $B\in
\co\tau_1$ and $A\cap B=\vnth$. Since $B\in \co\tau_1\sbs
\co\tau_2$ and $(Y,\tau_1',\tau_2')$ is $2\dd\WS$-supernormal in
$X$, there are neighborhoods $U'(A)\in \tau_2'$, $U(B)\in \tau_2$
such that $U'(A)\cap U(B)=\vnth$ and so $\tau_2\cl U'(A)\cap
U(B)=\vnth$. Furthermore, according to (2) of Corollary~2.3.12 in
[15],
$$  U'(A)\sbs \tau_2'\nt\tau_2'\cl U'(A)=\tau_1'\nt\tau_2'\cl U'(A)=
            U(A)\in \tau_1'     $$
and $U(A)\cap U(B)=\vnth$ so that $(Y,\tau_1',\tau_2')$ is
$(2,1)\dd\WS$-supernormal in $X$.

Finally, for the right vertical implication let $A\in \co\tau_1'$,
$B\in \co\tau_1$ and $A\cap B=\vnth$. Since $B\in \co\tau_1\sbs
\co\tau_2$ and $(Y,\tau_1',\tau_2')$ is $(1,2)\dd\WS$-supernormal
in $X$, there are neighborhoods $U'(A)\in \tau_2'$, $U(B)\in
\tau_1$ such that $U'(A)\cap U(B)=\vnth$. Hence $\tau_1\cl
U'(A)\cap U(B)=\vnth$ and by (2) of Corollary~2.3.12 in [15],
$$  U'(A)\sbs \tau_2'\nt\tau_2'\cl U'(A)=\tau_1'\nt\tau_2'\cl U'(A)=
            U(A)\in \tau_1'.     $$
Clearly $U(A)\cap U(B)=\vnth$ and, thus, $(Y,\tau_1',\tau_2')$ is
$1\dd\WS$-supernormal in $X$.

(2) First, let $A,\,B\in \co\tau_1'$ and $A\cap B=\vnth$. Since
$B\in\co\tau_1'\sbs \co\tau_2'$ and $(Y,\tau_1',\tau_2')$ is
$p\,$-strongly normal in $X$, there are neighborhoods $U'(A)\in
\tau_2$, $U(B)\in \tau_1$ such that $U'(A)\cap U(B)=\vnth$.
Clearly, $\tau_1\cl U'(A)\cap U(B)=\vnth$ and according to (2) of
Corollary~2.2.7 in [15],
$$  U'(A)\sbs \tau_2\nt\tau_1\cl U'(A)=\tau_1\nt\tau_1\cl U'(A)=
            U(A)\in \tau_1,     $$
where $U(A)\cap U(B)=\vnth$. Thus $(Y,\tau_1',\tau_2')$ is
$1$-strongly normal in $X$.

Furthermore, let $A\in \co\tau_1'$, $B\in \co\tau_2'$ and $A\cap
B=\vnth$. Since $A\in \co\tau_1'\sbs \co\tau_2'$ and
$(Y,\tau_1',\tau_2')$ is $2$-strongly normal in $X$, there are
neighborhoods $U(A),\,U'(B)\in \tau_2$ such that $U(A)\cap
U'(B)=\vnth$ and so $U(A)\cap\tau_2\cl U'(B)=\vnth$. Hence, by (2)
of Corollary~2.3.12 in [15],
$$  U'(B)\sbs \tau_2\nt\tau_2\cl U'(B)=\tau_1\nt\tau_2\cl U'(B)=
            U(B)\in \tau_1, \;\; U(A)\cap U(B)=\vnth     $$
and, hence, $(Y,\tau_1',\tau_2')$ is $p\,$-strongly normal in $X$.
The horizontal implications follow from the first implication of
(2), since $\tau_1<_N\tau_2$ implies $\tau_1<_C\tau_2$.~$\Box$
\vskip+0.2cm

Below we shall study the interrelation of the notions of relative
$\WS$-su\-per\-nor\-ma\-li\-ty and relative strong normality for
both the topological and the bitopological cases.

\begin{proposition}{2.19}
For a $\BsS$ $(Y,\tau_1',\tau_2')$ of a $\BS$ $(X,\tau_1,\tau_2)$
the following conditions are satisfied:
\begin{enumerate}
\item[$(1)$] If $Y\in \co\tau_i$ and $(Y,\tau_1',\tau_2')$ is
$p\,$-strongly normal in $X$, then
$$  (Y,\tau_1',\tau_2') \;\mbox{is $(i,j)\dd\WS$-supernormal in}\; X.  $$

\item[$(2)$] If $Y\in \tau_j$ and $(Y,\tau_1',\tau_2')$ is
$(i,j)\dd\WS$-supernormal in $X$, then
$$  (Y,\tau_1',\tau_2') \;\mbox{is $p\,$-strongly normal in}\; X. $$
\end{enumerate}
\end{proposition}

\noindent \textbf{Proof.\/} (1) Let $Y\in \co\tau_i$, $A\in
\co\tau_i'$, $B\in \co\tau_j$ and $A\cap B=\vnth$. Then $B'=B\cap
Y\in \co\tau_j'$ and since $(Y,\tau_1',\tau_2')$ is $p\,$-strongly
normal in $X$, there are disjoint sets $U'\in \tau_j$, $V'\in
\tau_i$ such that $A\sbs U'$ and $B'\sbs V'$. Evidently, if
$U=U'\cap Y$, $V=V'\cup (X\setminus Y)$, then $U\in \tau_j'$,
$V\in \tau_i$, $A\sbs U$, $B\sbs V$ and $U\cap V=\vnth$ so that
$(Y,\tau_1',\tau_2')$ is $(i,j)\dd\WS$-supernormal in $X$.

(2) Let $Y\in \tau_j$, $A\in \co\tau_i'$, $B\in\co\tau_j'$ and
$A\cap B=\vnth$. Then $A\cap \tau_j\cl B=\vnth$ and since
$(Y,\tau_1',\tau_2')$ is $(i,j)\dd\WS$-supernormal in $X$, there
are sets $U\in \tau_j'$, $V\in \tau_i$ such that $A\sbs U$,
$\tau_j\cl B\sbs V$ and $U\cap V=\vnth$. Clearly $U\in \tau_j$,
$B\sbs V$ and, hence, $(Y,\tau_1',\tau_2')$ is $p\,$-strongly
normal in $X$. $\Box$

\begin{corollary}{2.20}
For a $\BsS$ $(Y,\tau_1',\tau_2')$ of a $\BS$ $(X,\tau_1,\tau_2)$
the following conditions are satisfied:
\begin{enumerate}
\item[$(1)$] If $Y\in \co\tau_1\cap \co\tau_2$ and
$(Y,\tau_1',\tau_2')$ is $p\,$-strongly normal in $X$, then
$$  (Y,\tau_1',\tau_2') \;\mbox{is $p\dd\WS$-supernormal in}\; X. $$

\item[$(2)$] If $Y\in \tau_1\cap\tau_2$ and $(Y,\tau_1',\tau_2')$
is $p\,\dd\WS$-supernormal in $X$, then
$$  (Y,\tau_1',\tau_2') \;\mbox{is $p\,$-strongly normal in}\; X. $$

\item[$(3)$] If $Y\in \tau_1\cap\tau_2\cap\co\tau_1\cap\co\tau_2$,
then the following conditions are equivalent:
$(Y,\tau_1',\tau_2')$ is $(1,2)\dd\WS$-supernormal in $X$,
$(Y,\tau_1',\tau_2')$ is $(2,1)\dd\WS$-supernormal in $X$,
$(Y,\tau_1',\tau_2')$ is $p\,\dd\WS$-supernormal in $X$ and
$(Y,\tau_1',\tau_2')$ is $p\,$-strongly normal in $X$.
\end{enumerate}

Therefore, the equivalences remain valid for a $\BS$
$(X,\tau_1<\tau_2)$ if $Y\in \tau_1\cap\co\tau_1$.
\end{corollary}

\noindent \textbf{Proof.\/} (1) and (2) follow directly from (1)
and (2) of Proposition~2.19, respectively.

(3) By (1) and (2) of Proposition~2.19, if
$Y\in\tau_2\cap\co\tau_1$, then $(Y,\tau_1',\tau_2')$ is
$(1,2)\dd\WS$-supernormal in $X$ if and only if
$(Y,\tau_1',\tau_2')$ is $p\,$-strongly normal in $X$ and if
$Y\in\tau_1\cap\co\tau_2$, then $(Y,\tau_1',\tau_2')$ is
$(2,1)\dd\WS$-supernormal in $X$ if and only if
$(Y,\tau_1',\tau_2')$ is $p\,$-strongly normal in $X$.

The rest is obvious.~$\Box$

\begin{corollary}{2.21}
For a $\TsS$ $(Y,\tau')$ of a $\TS$ $(X,\tau)$ the following
conditions are satisfied:
\begin{enumerate}
\item[$(1)$] If $Y\in \co\tau$ and $(Y,\tau')$ is strongly normal
in $X$, then $(Y,\tau')$ is $\WS$-su\-per\-nor\-mal in $X$ and if
$Y\in \tau$ and $(Y,\tau')$ is $\WS$-supernormal in $X$, then
$(Y,\tau')$ is strongly normal in $X$.

\item[$(2)$] If $Y\in \tau\cap \co\tau$, then $(Y,\tau')$ is
$\WS$-supernormal in $X$ if and only if $(Y,\tau')$ is strongly
normal in $X$.
\end{enumerate}
\end{corollary}

In connection with Corollary~2.21 and so, with the bitopological
case, we consider it necessary to give the following elementary

\begin{example}{2.22}
\rm Let $X=\{a,b,c,d,e\}$,
$\tau=\{\vnth,\{c\},\{a,b,c\},\{c,d,e\},X\}$ and $Y\in \{a,b,d\}$.
Then $(Y,\tau')$ is $\WS$-supernormal in $X$, but it is not
strongly normal in $X$.

Now, if
$\tau=\{\vnth,\{a\},\{b\},\{a,b\},\{d,e\},\{a,d,e\},\{b,d,e\},\{a,b,d,e\},X\}$
and $Y=\{a,d,e\}$, then $(Y,\tau')$ is strongly normal in $X$, but
it is not $\WS$-supernormal in $X$.
\end{example}

Our next theorem, which is one of the main results of this work,
concerns the improvement of Theorem~3.2.26 and its
Corollaries~3.2.27--3.2.30 in [15].

\begin{theorem}{2.23}
If a $d$-second countable and $p\,$-normal $\BS$
$(X,\tau_1,\tau_2)$ can be represented as a union of a sequence
$F_1,F_2,\ldots$ of $p\,$-closed sets, where $F_k$ is \linebreak
$p\,\dd\WS$-supernormal in $X$ and $p\,\dd\ind F_k=0$ for each
$k=\ol{1,\infty}$, then $p\,\dd\Ind X=0$.
\end{theorem}

\noindent \textbf{Proof.\/} Let $A\in \co\tau_1$, $B\in \co\tau_2$
and $A\cap B=\vnth$. We shall prove that there exist sets $G\in
\tau_2$, $H\in \tau_1$ such that
\begin{equation}
    A\sbs G, \;\; B\sbs H, \;\; G\cap H=\vnth \;\;\mbox{and}\;\;
        G\cup H=X
\end{equation}
and so the empty set is a partition, corresponding to the pair
$(A,B)$.

Since $A\in \co\tau_1$, $B\in \co\tau_2$ and $A\cap B=\vnth$, by
Corollary~0.1.8 in [15], there exist sets $U_0\in \tau_2$, $V_0\in
\tau_1$ such that
\begin{equation}
    A\sbs U_0, \;\;\; B\sbs V_0 \;\; \mbox{and}\;\;
        \tau_1\cl U_0\cap \tau_2\cl V_0=\vnth.
\end{equation}

We shall define inductively two sequences of $2$-open and $1$-open
sets $U_0,U_1,\ldots$ and $V_0,V_1\ldots$, respectively,
satisfying for each $k=\ol{1,\infty}$ the following conditions:
\begin{eqnarray}
    &\ds U_{k-1}\sbs U_k, \;\; V_{k-1}\sbs V_k \;\;\mbox{if}\;\; k\geq 1,
        \;\;\mbox{and}\;\; \tau_1\cl U_k\cap \tau_2\cl V_k=\vnth \\
    &\ds F_k\sbs U_k\cup V_k, \;\;\mbox{where}\;\; F_0=\vnth.
\end{eqnarray}

Clearly, the sets $U_0$ and $V_0$, defined above, satisfy all
conditions for $k=0$. Assume that the sets $U_k$ and $V_k$,
satisfying (3) and (4), are defined for $k<p$. If $F_p=\tau_1\cl
F_p\cap \tau_2\cl F_p$, then the sets $\tau_1\cl U_{p-1}\cap
F_p\in\co\tau_1'$ and $\tau_2\cl V_{p-1}\cap F_p\in\co\tau_2'$ are
disjoint in $(F_p,\tau_1',\tau_2')$. Since $p\,\dd\ind F_p=0$,
according to Theorem~3.2.12 in [15], the empty set is a partition
between $\tau_1\cl U_{p-1}\cap F_p$ and $\tau_2\cl V_{p-1}\cap
F_p$ in $F_p$. Hence, there exists a set $V\in \tau_2'\cap
\co\tau_1'$ such that
$$  \tau_1\cl U_{p-1}\cap F_p\sbs V \;\; \mbox{and}\;\;
        \tau_2\cl V_{p-1}\cap F_p\sbs F_p\setminus V.   $$

It is clear that $\tau_1\cl V\cap \tau_2\cl (F_p\setminus
V)=\vnth$ as $\tau_1\cl V\sbs \tau_1\cl F_p$,
$\tau_2\cl(F_p\setminus V)\sbs \tau_2\cl F_p$ and $\tau_1\cl
F_p\cap \tau_2\cl F_p=F_p$. Moreover, $V\in \co\tau_1'$,
$\tau_2\cl V_{p-1}\in \co\tau_2$ and $V\cap \tau_2\cl
V_{p-1}=\vnth$. Since $(F_p,\tau_1',\tau_2')$ is
$(1,2)\dd\WS$-supernormal in $X$, there are disjoint sets $U\in
\tau_2'$, $W\in \tau_1$ such that $V\sbs U$ and $\tau_2\cl
V_{p-1}\sbs W$. Since $V\cap W=\vnth$ and $W\in \tau_1$, we have
$\tau_1\cl V\cap W=\vnth$ so that $\tau_1\cl V\cap \tau_2\cl
V_{p-1}\!=~\!\!\vnth$.

By the similar manner, taking into account that
$(F_p,\tau_1',\tau_2')$ is $(2,1)\dd\WS$-su\-per\-nor\-mal in $X$,
one can prove that $\tau_2\cl (F_p\setminus V)\cap \tau_1\cl
U_{p-1}=\vnth$.

Let $\Phi_1=\tau_1\cl U_{p-1}\cup \tau_1\cl V$ and
$\Phi_2=\tau_2\cl V_{p-1}\cup \tau_2\cl (F_p\setminus V)$. Then
$\Phi_i\in \co\tau_i$ and $\Phi_1\cap \Phi_2=\vnth$. Since
$(X,\tau_1,\tau_2)$ is $p\,$-normal, there exist sets $U_p\in
\tau_2$, $V_p\in \tau_1$ such that
\begin{eqnarray*}
    &\ds U_{p-1}\sbs \Phi_1\sbs U_p, \;\; V_{p-1}\sbs \Phi_2\sbs V_p, \\
    &\ds \tau_1\cl U_p\cap \tau_2\cl V_p=\vnth \;\;\mbox{and}\;\;
            F_p\sbs U_p\cup V_p
\end{eqnarray*}
so that the sets $\{U_p\}_{p=1}^\infty$ and $\{V_p\}_{p=1}^\infty$
satisfy (3) and (4) for $k=p$. Thus the construction of the
sequences $U_0,U_1,\ldots$ and $V_0,V_1,\ldots$ is completed. It
follows from (2), (3) and (4) that the unions
$G=\bigcup\limits_{p=1}^\infty U_p$ and
$H=\bigcup\limits_{p=1}^\infty V_p$ satisfy (1).~$\Box$

\begin{corollary}{2.24}
If a $p\,$-normal $\BS$ $(X,\tau_1,\tau_2)$ can be represented as
a union of a sequence $F_1,F_2,\ldots$ of $p\,$-closed sets, where
$F_k$ is $p\dd\WS$-supernormal in $X$ and $p\,\dd\Ind F_k=0$ for
each $k=\ol{1,\infty}$, then $p\,\dd\Ind X=0$.
\end{corollary}

The proof of this corollary is analogous to that of Theorem~2.23,
the only difference being that since $p\,\dd\Ind F_k=0$ for each
$k=\ol{1,\infty}$, Theorem~3.2.12 from [15] is unnecessary and so
is the requirement of $d$-second countability.

Note further that if $F_k\in\co\tau_1\cap\co\tau_2$, then the
requirement of $p\,\dd\WS$-su\-per\-nor\-ma\-li\-ty of each $F_k$
in $X$ is also unnecessary since it is automatically satisfied by
the reasoning before Example~2.17 provided that
$(X,\tau_1,\tau_2)$ is $p\,$-normal. Hence, in this case we obtain
Corollary~3.2.25 from [15].

\begin{corollary}{2.25}
If a $(d$-second countable$)$ $p\,$-normal $\BS$
$(X,\tau_1<\tau_2)$ can be represented as a union of a sequence
$F_1,F_2,\ldots$ of $2$-closed sets, where $F_k$ is
$(1,2)\dd\WS$-supernormal in $X$ and $(p\,\dd\ind F_k=0)$
$p\,\dd\Ind F_k=0$ for each $k=\ol{1,\infty}$, then $p\,\dd\Ind
X=0$.
\end{corollary}

\noindent \textbf{Proof.\/} On the one hand, $F_k\in \co\tau_2\sbs
p\dd\Cl(X)$ and on the other hand, $F_k\in \co\tau_2$ implies that
$F_k$ is $(2,1)\dd\WS$-supernormal in $X$. Hence, the case of
brackets follows directly from Theorem~2.23, and the case without
brackets -- from Corollary~2.24.~$\Box$

\begin{corollary}{2.26}
If a $(d$-second countable$)$ $p\,$-normal $\BS$
$(X,\tau_1,\tau_2)$ can be represented as a union of a sequence
$F_1,F_2,\ldots$, where every $F_k$ is a countable union of
$p\,$-closed sets, that is, $F_k=\bigcup\limits_{n=1}^\infty
F_n^k$, $(p\,\dd\ind F_k=0)$ $p\,\dd\Ind F_k=0$ for each
$k=\ol{1,\infty}$ and $F_n^k$ is $p\,\dd\WS$-supernormal in $X$
for each $n=\ol{1,\infty}$, $k=\ol{1,\infty}$, then $p\,\dd\Ind
X=0$. Moreover, if $(X,\tau_1,\tau_2)$ is $\R\dd\,p\,\dd\T_1$,
then $p\,\dd\ind X=0$.
\end{corollary}

\noindent \textbf{Proof.\/} For the case of brackets, note that by
(2) of Corollary~3.1.5 in [15], $p\,\dd\ind F_n^k=0$ and for the
case without brackets, by Corollary~3.2.8 in [15], $p\,\dd\Ind
F_n^k=0$ for each $n=\ol{1,\infty}$, $k=\ol{1,\infty}$. Therefore,
following Theorem~2.23 and Corollary~2.24, respectively,
$p\,\dd\Ind X=0$, since
$X=\bigcup\limits_{k=1}^\infty\bigcup\limits_{n=1}^\infty F_n^k$.

In the case where $(X,\tau_1,\tau_2)$ is $\R\dd\,p\,\dd\T_1$, it
suffices to use Corollary~3.2.13 from [15].~$\Box$

\begin{corollary}{2.27}
If a $\R\dd\,p\,\dd\T_1$, $d$-second countable and $p\,$-normal
$\BS$ \linebreak        $(X,\tau_1<\tau_2)$ can be represented as
a union of two $\BsS$'s $Y$ and $Z$, where $p\,\dd\Ind
Y=p\,\dd\Ind Z=0$ and one of them is $1$-open, then $p\,\dd\Ind
X=0$.
\end{corollary}

\noindent \textbf{Proof.\/} Let, for example, $Y\in \tau_1$. Then
$X\setminus Y\in \co\tau_1\sbs p\,\dd\Cl(X)$ and since $X\setminus
Y\sbs Z$, where $p\,\dd\Ind Z=0$, by Corollary~3.2.8 in [15],
$p\,\dd\Ind (X\setminus Y)=0$. Furthermore, by Corollary~0.1.13 in
[15], the $\BS$ $(X,\tau_1,\tau_2)$ is $p\,$-perfectly normal and,
hence, $Y\in \tau_1\sbs \tau_2\sbs 1\dd\cF_\sg(X)$, that is,
$Y=\bigcup\limits_{k=1}^\infty F_k$, where $F_k\in \co\tau_1\sbs
p\,\dd\Cl(X)$ for each $k=\ol{1,\infty}$. Thus, once more applying
Corollary 3.2.8 from [15], we obtain that $p\,\dd\Ind F_k=0$ for
each $k=\ol{1,\infty}$, since $p\,\dd\Ind Y=0$. Evidently,
$X=Y\cup (X\setminus Y)=\bigcup\limits_{k=1}^\infty F_k\cup
(X\setminus Y)$ and according to Corollary~2.24, $p\,\dd\Ind X=0$,
since $X\setminus Y\in\co\tau_1$ and $F_k\in\co\tau_1$ for each
$k=\ol{1,\infty}$ imply that each of them is $p\dd\WS$-supernormal
in $X$.~$\Box$

Below, in addition to Theorems 3.1.36 and 3.2.28 from [15], we
have

\begin{proposition}{2.28}
The following conditions are satisfied for a $\BS$
$(X,\tau_1,\tau_2)$:
\begin{enumerate}
\item[$(1)$] If $\tau_1C\tau_2$ and $(X,\tau_1,\tau_2)$ is
$(2,1)$-regular $($or $(2,1)\dd\ind(X)$ is finite$)$, then
$(2,1)\dd\ind X\leq 2\dd\ind X$.

\item[$(2)$] If $\tau_1N\tau_2$ and $(X,\tau_1,\tau_2)$ is
$2$-regular $($or $2\dd\ind (X)$ is finite$)$, then $1\dd\ind
X\leq (1,2)\dd\ind X$.

\item[$(3)$] If $\tau_1<_C\tau_2$ and $(X,\tau_1,\tau_2)$ is
$(2,1)$-regular or $\tau_1<_N\tau_2$ and $(X,\tau_1,\tau_2)$ is
$2$-regular, then
$$  1\dd\ind X=(1,2)\dd\ind X=2\dd\ind X=(2,1)\dd\ind X.    $$

\item[$(4)$] If $\tau_1<_C\tau_2$, $(X,\tau_1,\tau_2)$ is
$1\dd\T_1$ and $p$-normal $($or $(i,j)\dd\Ind X$ is finite$)$,
then
$$  1\dd\Ind X=(1,2)\dd\Ind X=2\dd\Ind X=(2,1)\dd\Ind X. $$
\end{enumerate}
\end{proposition}

\noindent \textbf{Proof.\/} (1) If $(X,\tau_1C\tau_2)$ is
$(2,1)$-regular, then by Corollary~2.2.9 in [15],
$\tau_2\sbs\tau_1$ and it remains to use the first inequality in
(1) of Theorem 3.1.36 in [15].

(2) If $(X,\tau_1N\tau_2)$ is $2$-regular, then by Theorem~1 in
[22], $\tau_2\sbs \tau_1$ and it remains to use the second
inequality in (1) of Theorem~3.1.36 in [15].

(3) If $(X,\tau_1<_C\tau_2)$ is $2$-regular, then by (1) or by (2)
above, $\tau_1=\tau_2$.

(4) If $(X,\tau_1<_C\tau_2)$ is $1\dd\T_1$ and $p$-normal, then by
the implications before Definition~0.1.7 in [15],
$(X,\tau_1,\tau_2)$ is $(2,1)$-regular and so, by (1) above,
$\tau_1=\tau_2$.

Finally, note that if $(2,1)\dd\ind X$ $(2\dd\ind X$,
$(i,j)\dd\Ind X)$ is finite, then by (1) of Proposition~3.1.4 in
[15] (the well-know topological fact, (1) of Proposition~3.2.7 in
[15]), the $\BS$ $(X,\tau_1,\tau_2)$ is $(2,1)$-regular
$(2$-regular, $p$-normal).~$\Box$ \vskip+0.2cm

At the end of this section we consider the notion of almost
$n$-dimensionality [2] from the bitopological point of view.

\begin{definition}{2.29}
\rm We say that a $\TS$ $(X,\tau_2)$ is almost $n$-dimensional in
the sense of the small (large) inductive dimension and we write
$\aind(X,\tau_2)=n$ $(\aInd(X,\tau_2)=n)$ if there exists a
topology $\tau_1$ on $X$such that $\ind(X,\tau_1)\leq n$
$(\Ind(X,\tau_1)\leq n)$, $\tau_1$ is weaker than the given
topology $\tau_2$ on $X$, the $\BS$ $(X,\tau_1<\tau_2)$ is
(2,1)-regular $(p$-normal) and $n$ is the smallest natural number
such that a topology $\tau_1$ exists for~$n$.

In this case we also say that the $\TS$ $(X,\tau_2)$ is almost
$n$-dimensional (in both senses) owing to the (weaker) topology
$\tau_1$ on $X$.
\end{definition}

\begin{theorem}{2.30}
For a $\TS$ $(X,\tau_2)$  we have:
\begin{enumerate}
\item[$(1)$] If $\aind(X,\tau_2)=n$ owing to $\tau_1$, then
$(1,2)\dd\ind(X,\tau_1,\tau_2)\leq n$, the $\BS$
$(X,\tau_1,\tau_2)$ is $d$-regular and $p$-regular, and $\tau_1$
is a cotopology of $\tau_2$ in the sense of $[1]$.

\item[$(2)$] If $\aInd (X,\tau_2)=n$ owing to $\tau_1$, then
$(X,\tau_1,\tau_2)$ is $1$-normal; moreover, if
$(X,\tau_1,\tau_2)$ is $1\dd\T_1$, then $\aind(X,\tau_2)\leq n$.

\item[$(3)$] If $\aind (X,\tau_2)\!=\!n$ owing to $\tau_1$, then
for each topology $\tau$ on $X$ such that
$\tau_1\!\sbs\!\tau\!\sbs~\!\!\tau_2$, the $\BS$ $(X,\tau,\tau_2)$
is $(2,1)$-regular, $2$-regular and
$(1,2)\dd\ind(X,\tau_1,\tau)\!\!\leq~\!\!\!n$, the $\BS$
$(X,\tau_1,\tau)$ is $1$-regular and $(1,2)$-regular.

\item[$(4)$] If $\aind (X,\tau_2)=n$ owing to $\tau_1$,
$(X,\tau_1,\tau_2)$ is $1\dd\T_2$ and $2$-locally compact, then
$(X,\tau_1,\tau_2)$ is a $2$-Baire space $($and, hence, it is also
an almost $(2,1)$-Baire space, a $2$-weak Baire space and a
$(2,1)$-weak Baire space in the sense of $[15])$; moreover, if, in
addition, $(X,\tau_1,\tau_2)$ is $1$-second countable, then the
$\TS$ $(X,\tau_2)$ has a bounded complete computational model in
the sense of $[8]$.
\end{enumerate}
\end{theorem}

\noindent \textbf{Proof.\/} (1) By Definition~2.29 and (1) of
Theorem~3.1.36 in [15],
$$  (1,2)\dd\ind(X,\tau_1,\tau_2)\leq
        1\dd\ind(X,\tau_1,\tau_2)\leq n,        $$
where the right inequality gives that $(X,\tau_1,\tau_2)$ is
$1$-regular. But it is also $(2,1)$-regular and, hence, it is
$d$-regular and $p$-regular as $\tau_1\sbs\tau_2$. The fact that
$\tau_1$ is a cotopology of $\tau_2$ follows directly from
Definition~2.29 and Theorem~7.3.20 in [15].

(2) Since $1\dd\Ind(X,\tau_1,\tau_2)\leq n$, $(X,\tau_1,\tau_2)$
is $1$-normal. If $(X,\tau_1,\tau_2)$ is $1\dd\T_1(\llra
\R\dd\,p\,\dd\T_1)$, then $(X,\tau_1,\tau_2)$ is $p$-normal
implies that $(X,\tau_1,\tau_2)$ is $p$-regular and so,
$(2,1)$-regular. Moreover,
$$  1\dd\ind(X,\tau_1,\tau_2)\leq
            1\dd\Ind(X,\tau_1,\tau_2)\leq n,    $$
i.e., $\ind(X,\tau_1)\leq \Ind(X,\tau_1)$ and by Definition~2.29,
$\aind(X,\tau_2)\leq n$.

(3) Follows directly from (1) above and (1) of Theorem~3.1.36 in
[15] as $\tau_1\sbs\tau_2$.

(4) The first part follows directly from (1) and (2) of
Corollary~7.3.25 in [15]. The second part is an immediate
consequence of (1) and (2) of Corollary~7.3.25 in [15] taking into
account (3) of Theorem~9 in [8], since, by (1) above,
$(X,\tau_1,\tau_2)$ is $p$-regular.~$\Box$ \vskip+0.2cm

\section*{3. $(i,j)$-Submaximal, $(i,j)$-Nodec, $(i,j)\dd\,\cI$
and} \section*{\ \ \ \ $(i,j)\dd\,D$-Spaces}

\begin{definition}{3.1}
\rm A $\BS$ $(X,\tau_1,\tau_2)$ is said to be $(i,j)$-submaximal
if every subset of $X$ is $(i,j)$-locally closed, that is, if
$2^X=(i,j)\dd\LC(X)$.
\end{definition}

Taking into account the inclusions after Definition~1.3, for a
$\BS$ $(X,\tau_1<\tau_2)$ the following implications hold:
$$  \begin{matrix}
        (X,\tau_1,\tau_2)\;\mbox{is $1$-submaximal} & \lra &
            (X,\tau_1,\tau_2)\;\mbox{is $(2,1)$-submaximal} \\
        \big\Downarrow & & \big\Downarrow \\
        (X,\tau_1,\tau_2)\;\mbox{is $(1,2)$-submaximal} & \lra &
            (X,\tau_1,\tau_2)\;\mbox{is $2$-submaximal}\,.
    \end{matrix}        $$

Therefore, from Corollary 1.3 in [5] it follows immediately that
if \linebreak        $(X,\tau_1<\tau_2)$ is $i$-submaximal or
$(i,j)$-submaximal, then it is $2$-nodec.

\begin{theorem}{3.2}
For a $\BS$ $(X,\tau_1,\tau_2)$ the following conditions are
equivalent:
\begin{enumerate}
\item[$(1)$] $(X,\tau_1,\tau_2)$ is $(i,j)$-submaximal.

\item[$(2)$] Every subset of $X$ is $\co(i,j)$-locally closed.

\item[$(3)$] Every $j$-boundary subset of $X$ is $i$-closed.

\item[$(4)$] $\tau_j\cl A\setminus A$ is $i$-closed for every
subset $A\sbs X$.

\item[$(5)$] Every $j$-dense subset of $X$ is $i$-open.
\end{enumerate}
\end{theorem}

\noindent \textbf{Proof.\/} It is clear that $(1)\llra (2)$ since
$2^X=(i,j)\dd\LC(X)$ implies the equivalence $A\in
(i,j)\dd\LC(X)\llra X\setminus A\in (i,j)\dd\LC(X)$ so that $A\in
\co(i,j)\dd\LC(X)$ for each subset $A\sbs X$.

$(4)\llra (5)$ If $A\in j\dd\,\cD(X)$, then, by (4), $\tau_j\cl
A\setminus A=X\setminus A\in \co\tau_i$ so that $A\in \tau_i$.
Conversely, if there is a set $A\sbs X$ such that $(\tau_j\cl
A\setminus A)\,\ol{\in}\,\co\tau_i$, then
$$  X\setminus (\tau_j\cl A\setminus A)=
        \tau_j\nt(X\setminus A)\cup A\,\ol{\in}\,\tau_i.    $$
Hence, by (5), $\tau_j\nt(X\setminus A)\cup
A\,\ol{\in}\,j\dd\,\cD(X)$.

$(4)\lra (1)$ Let $A\sbs X$ be any set. Then, by (4), $\tau_j\cl
A\setminus A=\tau_j\cl A\cap (X\setminus A)\in \co\tau_i$ so that
$$  X\setminus (\tau_j\cl A\setminus A)=
        (X\setminus \tau_j\cl A)\cup A\in \tau_i.       $$
But
$$  A=\tau_j\cl A\cap((X\setminus \tau_j\cl A)\cup A)=F\cap U,  $$
where $F\in \co\tau_j$ and $U\in \tau_i$; hence $A\in
(i,j)\dd\LC(X)$.

$(2)\lra (3)$ If $\tau_j\nt A=\vnth$, then the equivalence $A\in
\co(i,j)\dd\LC(X)\llra (A=U\cup F$, where $U\in \tau_j$ and $F\in
\co\tau_i)$ implies that $U=\vnth$ and, thus, $A\in \co\tau_i$.

The implication $(3)\lra (4)$ is obvious, since $\tau_j\cl
A\setminus A\in j\dd\Bd(X)$ for each subset $A\sbs X$.~$\Box$
\vskip+0.2cm

Note also here that according to (3) of Theorem~3.2 above and
Corollary~1.3 in [5], for a $\BS$ $(X,\tau_1S\tau_2)$ we have:
\begin{eqnarray*}
    &\ds (X,\tau_1,\tau_2)\;\mbox{is $(1,2)$-submaximal}\llra
        (X,\tau_1,\tau_2)\;\mbox{is $1$-submaximal}\lra \\
    &\ds \lra (X,\tau_1,\tau_2)\;\mbox{is $1$-nodec}
\end{eqnarray*}
and
\begin{eqnarray*}
    &\ds (X,\tau_1,\tau_2)\;\mbox{is $(2,1)$-submaximal}\llra
        (X,\tau_1,\tau_2)\;\mbox{is $2$-submaximal}\lra \\
    &\ds \lra (X,\tau_1,\tau_2)\;\mbox{is $2$-nodec}\,,
\end{eqnarray*}
since, by (2) of Theorem~2.1.5 in [15], if $\tau_1S\tau_2$, then
$1\dd\Bd(X)=2\dd\Bd(X)$.

The other relations between submaximal and nodec spaces as well as
between their bitopological modifications will be given in
Corollary~3.7.

\begin{corollary}{3.3}
For a $j\dd\T_1$ $\BS$ $(X,\tau_1,\tau_2)$ the following
conditions are equivalent:
\begin{enumerate}
\item[$(1)$] Every $j$-boundary subset of $X$ is $i$-discrete.

\item[$(2)$] $\tau_j\cl A\setminus A$ is $i$-discrete for every
subset $A\sbs X$, \\
and in the case where $(X,\tau_1,\tau_2)$ is $(i,j)$-submaximal,
each of them is satisfied.
\end{enumerate}

Moreover, if $X_j^i=\vnth$, then $(1)$ $($and so $(2))$ is
equivalent to the $(i,j)$-sub\-ma\-xi\-ma\-li\-ty of
$(X,\tau_1,\tau_2)$.
\end{corollary}

\noindent \textbf{Proof.\/} It is clear that for any subset $A\sbs
X$ the set $\tau_j\cl A\setminus A\in j\dd\Bd(X)$ and, hence, by
(1), $\tau_j\cl A\setminus A=(\tau_j\cl A\setminus A)_i^i$.
Conversely, if $A\in j\dd\Bd(X)$ (i.e., if $X\setminus A\in
j\dd\,\cD(X))$, then by (2),
$$  \tau_j\cl (X\setminus A)\setminus (X\setminus A)=
        \big(\tau_j\cl (X\setminus A)\setminus
                (X\setminus A)\big)_i^i,     $$
that is, $X\setminus (X\setminus A)=A=A_i^i$.

Now, let $(X,\tau_1,\tau_2)$ be $(i,j)$-submaximal and let us
prove that (3) of Theorem~3.2 implies (1) of this corollary.
Indeed: if there is a set $A\in j\dd\Bd(X)$ such that $A\neq
A_i^i$, i.e., $A\setminus A_i^i\neq\vnth$, then there exists a
point $x\in A$ such that for each $i$-open neighborhood $U(x)$:
$U(x)\cap (A\setminus \{x\})\neq \vnth$. If $B=A\setminus \{x\}$,
then $x\in \tau_i\cl B$ so that $B\neq \tau_i\cl B$. But $B\sbs A$
implies that $B\in j\dd\Bd(X)$. A contradiction with (3) of
Theorem~3.2.

Finally, let us prove that if $X_j^i=\vnth$, then (1) of this
corollary implies (3) of Theorem~3.2.

Let $A\sbs X$, $A\in j\dd\Bd(X)$ and $A\neq \tau_i\cl A$. If $x\in
\tau_i\cl A\setminus A$ is any point, then for each $i$-open
neighborhood $U(x)$ we have $U(x)\cap (A\setminus \{x\})=U(x)\cap
A\neq \vnth$. Let $B=A\cup \{x\}$, where $\{x\}\in \co\tau_j$ as
$(X,\tau_1,\tau_2)$ is $j\dd\T_1$. Then, by the well-known
topological fact,
$$  \tau_j\nt B=\tau_j\nt(A\cup \{x\})=
        \tau_j\nt\big(\tau_j\nt A\cup\{x\}\big)=
                \tau_j\nt\{x\}=\vnth     $$
as $X_j^i=\vnth$ so that, by (1), $B=B_i^i$, which is impossible,
since for each $i$-open neighborhood $U(x)$: $U(x)\cap (B\setminus
\{x\})=U(x)\cap A\neq\vnth$.~$\Box$ \vskip+0.2cm

Note also here, that if $X_j^i\neq\vnth$, then (1) (and so (2)) of
Corollary 3.3 does not imply the $(i,j)$-submaximality of $X$.

\begin{example}{3.4}
\rm Let $X=\{a,b,c\}$, $\tau_1=\{\vnth,\{b\},\{c\},\{b,c\},X\}$
and $\tau_2=        \linebreak
   \{\vnth,\{a\},\{a,b\},\{a,c\},X\}$. It is clear that
$2\dd\Bd(X)=\{\{b\},\{c\},\{b,c\}\}$ and each of these sets is
$1$-discrete, while the set $A=\{a,b\}\,\ol{\in}\,(1,2)\dd\LC(X)$
so that $(X,\tau_1,\tau_2)$ is not $(1,2)$-submaximal. The reason
of this is that $X_2^i=\{\{a\}\}\neq\vnth$.
\end{example}

\begin{corollary}{3.5}
For a $1\dd\T_1$ $\BS$ $(X,\tau_1<\tau_2)$ the following
conditions are equivalent:
\begin{enumerate}
\item[$(1)$] $(X,\tau_1,\tau_2)$ is $(2,1)$-submaximal.

\item[$(2)$] Every subset of $X$ is $\co(2,1)$-locally closed.

\item[$(3)$] If $A\in 1\dd\Bd(X)$, then $A\in \co\tau_2$.

\item[$(4)$] If $A\in 1\dd\Bd(X)$, then $A=A_2^i$.

\item[$(5)$] $\tau_1\cl A\setminus A\in\co\tau_2$ for every subset
$A\sbs X$.

\item[$(6)$] $\tau_1\cl A\setminus A=(\tau_1\cl A\setminus A)_2^i$
for every subset $A\sbs X$.

\item[$(7)$] If $A\in 1\dd\,\cD(X)$, then $A\in \tau_2$.
\end{enumerate}
\end{corollary}

\noindent \textbf{Proof.\/} Taking into account Theorem~3.2 and
the first part of Corollary~3.3, it suffices to prove only that
$(6)\lra (5)$. Let $A\sbs X$ be a set such that \linebreak
   $\tau_1\cl A\setminus A\,\ol{\in}\,\co\tau_2$ so that there exists
a point $p\in (\tau_1\cl A\setminus A)_2^d\setminus (\tau_1\cl
A\setminus A)$. Then
$$  p\in (\tau_1\cl A\setminus A)_2^d\sbs
        \tau_2\cl(\tau_1\cl A\setminus A)\sbs
            \tau_2\cl \tau_1\cl A=\tau_1\cl A.      $$
Hence $p\,\ol{\in}\,X\setminus \tau_1\cl A$. On the other hand,
$p\,\ol{\in}\,(\tau_1\cl A\setminus A)$ implies that
  \linebreak     $p\in (X\setminus \tau_1\cl A)\cup A$ and, thus, $p\in A$. Let
$B=A\setminus \{p\}$. Then $\tau_1\cl A=\tau_1\cl B\cup\{p\}$ as
$(X,\tau_1,\tau_2)$ is $1\dd\T_1(\llra \R\dd\,p\dd\,\T_1)$.
Moreover, $p\in (\tau_1\cl A\setminus A)_2^d\sbs (\tau_1\cl
A)_2^d\sbs (\tau_1\cl A)_1^d=A_1^d$ so that for any $1$-open
neighborhood $U(p):\;U(p)\cap (A\setminus \{p\})=U(p)\cap
B\neq\vnth$. Hence $p\in \tau_1\cl B$ and, thus, $\tau_1\cl
B=\tau_1\cl A$. Finally,
$$  \tau_1\cl B\setminus B=\tau_1\cl A\setminus (A\setminus \{p\})=
        (\tau_1\cl A\setminus A)\cup \{p\}.     $$
If there is a $2$-open neighborhood $U(p)$ such that
$$  U(p)\cap (\tau_1\cl B\setminus B)=U(p)\cap
        \big((\tau_1\cl A\setminus A)\cup \{p\}\big)=\{p\}, $$
then $U(p)\cap (\tau_1\cl A\setminus A)=\vnth$ which contradicts
the assumption that    \linebreak     $p\in \tau_2\cl(\tau_1\cl
A\setminus A)$.~$\Box$

\begin{corollary}{3.6}
Every $\BsS$ $(Y,\tau_1',\tau_2')$ of an $(i,j)$-submaximal $\BS$
$(X,\tau_1,\tau_2)$ is also $(i,j)$-submaximal and if $Y$ is
$j$-discrete and $Y\cap X_j^i=\vnth$, then $Y$ is $i$-closed.
\end{corollary}

\noindent \textbf{Proof.\/} The first part is obvious. For the
second part, by (3) of Theorem~3.2, it suffices to prove only that
$\tau_j\nt Y=\vnth$. Contrary: $\tau_j\nt Y\neq\vnth$; then there
are a point $x\in Y$ and a $j$-open neighborhood $U(x)$ such that
$U(x)\sbs Y$. Since $Y=Y_j^i$, there is a $j$-open neighborhood
$V(x)$ such that $V(x)\cap Y=\{x\}$. Hence, $V(x)\cap U(x)=\{x\}$
which contradicts the condition $Y\cap X_j^i=\vnth$.~$\Box$

\begin{corollary}{3.7}
For a $\BS$ $(X,\tau_1<\tau_2)$ we have:
\begin{enumerate}
\item[$(1)$] If $(X,\tau_1,\tau_2)$ is $1$-submaximal, then it is
$d$-nodec and $p\,$-nodec.

\item[$(2)$] If $(X,\tau_1,\tau_2)$ is $(1,2)$-submaximal, then it
is $(2,1)$-nodec.
\end{enumerate}

For a $\BS$ $(X,\tau_1<_C\tau_2)$ we have
\begin{enumerate}
\item[$(3)$] If $(X,\tau_1,\tau_2)$ is $(1,2)$-submaximal, then it
is $1$-nodec.
\end{enumerate}

For a $\BS$ $(X,\tau_1<_N\tau_2)$ as well as for a $\BS$
$(X,\tau_1<_S\tau_2)$ we have
\begin{enumerate}
\item[$(4)$] If $(X,\tau_1,\tau_2)$ is $(1,2)$-submaximal, then it
is $d$-nodec and $p\,$-nodec.
\end{enumerate}
\end{corollary}

\noindent \textbf{Proof.\/} (1) Let $A\in i\dd\ND(X)\cup
(i,j)\dd\ND(X)$. Then, in all cases, $A\in 1\dd\Bd(X)$ as
$\tau_1\sbs \tau_2$. Hence, by (c) and (d) of Theorem~1.2 in [5],
$A=\tau_1\cl A=A_1^i$. But $\tau_1\sbs\tau_2$ implies that
$A_1^i\sbs A_2^i\sbs A\sbs \tau_2\cl A\sbs\tau_1\cl A$ and so
$A_1^i=A_2^i=A=\tau_2\cl A=\tau_1\cl A$.

Therefore, $(X,\tau_1,\tau_2)$ is $d$-nodec and $p\,$-nodec.

(2) If $A\in (2,1)\dd\ND(X)$, then $A\in 2\dd\Bd(X)$ and by (3) of
Theorem~3.2, $A=\tau_1\cl A$. Moreover, by the first part of
Corollary~3.3, $A=A_1^i$ and, hence, $A=A_1^i=A_2^i=\tau_1\cl A$
so that $A=\tau_1\cl A=A_2^i$. Thus $(X,\tau_1,\tau_2)$ is
$(2,1)$-nodec.

(3) Let $A\!\in\!1\dd\ND(X)$. Then, by (3) of Theorem~2.2.20 in
[15], $A\!\in\!(2,1)\dd\ND(X)$ and by the same reasonings as in
the proof of (2), we obtain that $A=\tau_1\cl A=A_1^i$, that is,
$(X,\tau_1,\tau_2)$ is $1$-nodec.

(4) For the case where $\tau_1<_S\tau_2$, we have
$(X,\tau_1,\tau_2)$ is $(1,2)$-submaximal $\llra
(X,\tau_1,\tau_2)$ is $1$-submaximal and, thus, it remains to use
(1).

Let $\tau_1<_N\tau_2$ and $A\in i\dd\ND(X)\cup (i,j)\dd\ND(X)$.
Clearly, in all cases, $A\in (1,2)\dd\ND(X)$ and according to (3)
of Theorem~2.3.19 in [15], $A\in 2\dd\ND(X)$. Therefore, $A\in
2\dd\Bd(X)$ and since $(X,\tau_1,\tau_2)$ is $(1,2)$-submaximal,
by (3) of Theorem~3.2 and the first part of Corollary~3.3,
$A=\tau_1\cl A=A_1^i$. The rest is obvious.~$\Box$

\begin{definition}{3.8}
\rm A $\BS$ $(X,\tau_1,\tau_2)$ is $(i,j)$-strongly $\sg$-discrete
if it can be represented as the union of a countable family of
$j$-closed $i$-discrete subspaces, that is, if
$X=\bigcup\limits_{n=1}^\infty A_n$, where $A_n=\tau_j\cl
A_n=(A_n)_i^i$ for each $n=\ol{1,\infty}$.
\end{definition}

Clearly, for a $\BS$ $(X,\tau_1<\tau_2)$ the following
implications hold:
$$  \begin{matrix}
        (X,\tau_1,\tau_2)\;\mbox{is $1$-strongly $\sg$-discrete}
            \!&\!\!\lra\!\!&\!(X,\tau_1,\tau_2)\;\mbox{is $(2,1)$-strongly $\sg$-discrete} \\
        \big\Downarrow & & \big\Downarrow \\
        (X,\tau_1,\tau_2)\;\mbox{is $(1,2)$-strongly $\sg$-discrete}
            \!&\!\!\lra\!\!&\!(X,\tau_1,\tau_2)\;\mbox{is $2$-strongly $\sg$-discrete}\,.
    \end{matrix}        $$

Now, it is not difficult to see that take place

\begin{theorem}{3.9}
Every $(1,2)$-nodec $\BS$ $(X,\tau_1<\tau_2)$ for which  no
non\-empty $1$-open subset is $(1,2)$-strongly $\sg$-discrete, is
a $(1,2)\dd\BrS$, and hence, a $1\dd\BrS$.
\end{theorem}

\noindent \textbf{Proof.\/} Let $U\in \tau_1\setminus \{\vnth\}$
and $U$ is of $(1,2)\dd\Catg\I$. Then, by Corollary~1.5.14 in
[15], $U\in (1,2)\dd\Catg_{{}_{\I}}(X)$ so that
$U=\bigcup\limits_{n=1}^\infty A_n$, where $A_n\in (1,2)\dd\ND(X)$
for each $n=\ol{1,\infty}$. Evidently, each set $A_n$ is
$2$-closed and $1$-discrete, since $(X,\tau_1,\tau_2)$ is
$(1,2)$-nodec and, thus, $U$ is $(1,2)$-strongly $\sg$-discrete.
The rest follows from (1) of Theorem~4.1.3 in [15].~$\Box$

\begin{corollary}{3.10}
For a $\BS$ $(X,\tau_1<_C\tau_2)$ the following conditions are
sa\-tis\-fi\-ed:
\begin{enumerate}
\item[$(1)$] If $(X,\tau_1,\tau_2)$ is $(1,2)$-nodec and a set
$U\in \tau_1\setminus \{\vnth\}$ is $(1,2)$-strongly
$\sg$-discrete, then $\tau_2\cl U\;(\llra \tau_1\cl U)$ is
$2$-strongly $\sg$-discrete.

\item[$(2)$] If $(X,\tau_1,\tau_2)$ is a $(2,1)$-nodec space for
which no nonempty $1$-open $(2$-open$)$ subset is $(2,1)$-strongly
$\sg$-discrete, then it is a $1\dd\BrS$ $($an $A\dd(2,1)\dd\BrS)$.
\end{enumerate}
\end{corollary}

\noindent \textbf{Proof.\/} (1) Evidently, $(X,\tau_1,\tau_2)$ is
$2$-nodec and if $U\in \tau_1\setminus \{\vnth\}$, $U$ is
\linebreak         $(1,2)$-strongly $\sg$-discrete, then $U$ is
$2$-strongly $\sg$-discrete. Hence,  $\tau_2\cl U$ is also
$2$-strongly $\sg$-discrete, since $U\in \tau_1\sbs\tau_2$. On the
other hand, according to (3) of Corollary~2.2.7 in [15],
$\tau_2\cl U=\tau_1\cl U$.

(2) Contrary: there is a set $U\in \tau_1\setminus \{\vnth\}$ such
that $U$ is of $1\dd\Catg\I\;(\llra U\in 1\dd\Catg_{{}_{\I}}(X))$.
Then $U=\bigcup\limits_{n=1}^\infty A_n$, where $A_n\in
1\dd\ND(X)$ and, hence, by (3) of Theorem~2.2.20 in [15], $A_n\in
(2,1)\dd\ND(X)$ for each $n=\ol{1,\infty}$. Therefore,
$(X,\tau_1,\tau_2)$ is $(2,1)$-nodec implies that $A_n=\tau_1\cl
A_n=(A_n)_2^i$ for each $n=\ol{1,\infty}$ and thus, $U$ is
$(2,1)$-strongly $\sg$-discrete.

Now, taking into account (6) of Definition~1.1, the case of
brackets is obvious.~$\Box$ \vskip+0.2cm

Clearly, (2) of Corollary~3.10 and also Proposition~4.14 in [5]
contain the sufficient conditions for a $\BS$ to be $1$-Baire.

Very close to (1) of Corollary~3.10 is the following

\begin{proposition}{3.11}
If a $p\,$-normal $\BS$ $(X,\tau_1<_C\tau_2)$ is $(1,2)$-nodec and
a set $U\in \tau_1\setminus \{\vnth\}$ is $(1,2)$-strongly
$\sg$-discrete, then  $\tau_2\cl U\;(\llra\tau_1\cl U)$ is also
\linebreak       $(1,2)$-strongly $\sg$-discrete.
\end{proposition}

\noindent \textbf{Proof.\/} Let $U\in \tau_1\setminus \{\vnth\}$
and $U=\bigcup\limits_{n=1}^\infty A_n$, where $A_n=\tau_2\cl
A_n=(A_n)_1^i$ for each $n=\ol{1,\infty}$. It is clear that
$\tau_2\cl U\setminus U\in (1,2)\dd\ND(X)$.

Since $(X,\tau_1,\tau_2)$ is $(1,2)$-nodec,
$$  \tau_2\cl U\setminus U=\tau_2\cl(\tau_2\cl U\setminus U)=
        (\tau_2\cl U\setminus U)_1^i        $$
so that
$$  \tau_2\cl U=U\cup(\tau_2\cl U\setminus U)=
        {\textstyle\bigcup\limits_{n=1}^\infty} B_n, $$
where $B_n=A_n\cup (\tau_2\cl U\setminus U)$ and
$$  \tau_2\cl B_n=\tau_2\cl A_n\cup
            \tau_2\cl(\tau_2\cl U\setminus U)=
        A_n\cup (\tau_2\cl U\setminus U)=B_n     $$
for each $n=\ol{1,\infty}$. Hence, it remains to prove only that
$B_n=(B_n)_1^i$, that is,
$$  A_n\cup (\tau_2\cl U\setminus U)=
        \big(A_n\cup (\tau_2\cl U\setminus U)\big)_1^i     $$
for each $n=\ol{1,\infty}$. Since $A_n=(A_n)_1^i$ and $\tau_2\cl
U\setminus U=(\tau_2\cl U\setminus U)_1^i$, it suffices to prove
only that
$$  \big(A_n\cup (\tau_2\cl U\setminus U)\big)_1^i=(B_n)_1^i=
        (A_n)_1^i\cup (\tau_2\cl U\setminus U)_1^i  $$
for each $n=\ol{1,\infty}$. The equivality $A_n\cap (\tau_2\cl
U\setminus U)=\vnth$ implies that $(B_n)_1^i\sbs (A_n)_1^i\cup
(\tau_2\cl U\setminus U)_1^i$ for each $n=\ol{1,\infty}$. Indeed:
for arbitrary $n\in N$ let $x\in (B_n)_1^i$; then there is a
$1$-open neighborhood $U(x)$ such that $U(x)\cap B_n=\{x\}$. Since
$B_n=A_n\cup (\tau_2\cl U\setminus U)$, where $A_n\cap (\tau_2\cl
U\setminus U)=\vnth$, we have $U(x)\cap A_n=\{x\}$ or $U(x)\cap
(\tau_2\cl U\setminus U)=\{x\}$ so that $x\in (A_n)_1^i$ or $x\in
(\tau_2\cl U\setminus U)_1^i$.

Finally, let us prove that $(A_n)_1^i\cup (\tau_2\cl U\setminus
U)_1^i\sbs (B_n)_1^i$. First, let $x\in (A_n)_1^i$. Then there is
a $1$-open neighborhood $U(x)$ such that $U(x)\cap A_n=\{x\}$.
Since          \linebreak           $\tau_1<_C\tau_2$ implies that
$\tau_1\cl U=\tau_2\cl U$, we have that $V(x)=U(x)\setminus
(\tau_2\cl U\setminus U)\in \tau_1$, where $V(x)\cap B_n=V(x)\cap
A_n=\{x\}$ so that $x\in (B_n)_1^i$. Now, let $x\in (\tau_2\cl
U\setminus U)_1^i$. Then there is a neighborhood $W(x)\in \tau_1$
such that $W(x)\cap (\tau_2\cl U\setminus U)=\{x\}$. Hence
$x\,\ol{\in}\,U=\bigcup\limits_{n=1}^\infty A_n$ so that
$x\,\ol{\in}\,A_n$ for each $n=\ol{1,\infty}$. Since
$(X,\tau_1,\tau_2)$ is       \linebreak  $p\,$-normal, $A_n\in
\co\tau_2$, $U\in \tau_1$ and $A_n\sbs U$, there is a set $V_n\in
\tau_1$ such that $A_n\sbs V_n\sbs \tau_2\cl V_n=\tau_1\cl V_n\sbs
U$, and $x\,\ol{\in}\,\tau_1\cl V_n$ for each $n=\ol{1,\infty}$.
Putting now $V_n(x)=W(x)\setminus \tau_1\cl V_n$ we obtain that
$V_n(x)\in \tau_1$ and
\begin{eqnarray*}
       V_n(x)\cap B_n &&\hskip-0.6cm =(V_n(x)\cap A_n)\cup \big(V_n(x)\cap
            (\tau_2\cl U\setminus U)\big)= \\
      &&\hskip-0.6cm =V_n(x)\cap (\tau_2\cl U\setminus U)=\{x\}
\end{eqnarray*}
so that $x\in (B_n)_1^i$ for each $n=\ol{1,\infty}$.~$\Box$

\begin{proposition}{3.12}
For a $(2,1)$-submaximal $\BS$ $(X,\tau_1<\tau_2)$ the following
conditions are satisfied:
\begin{enumerate}
\item[$(1)$] If $(Y,\tau_1',\tau_2')$ is a $\BsS$ of $X$ and
$U=\tau_1\nt Y$, then the set $Z=Y\setminus Y_2^i\sbs\tau_2\cl U$
and $Y=A\cup B$, where $A\in \tau_1$ and $B=\tau_2\cl B=B_2^i$.

\item[$(2)$] If $(C,\tau_1',\tau_2')$ is a $2$-connected $\BsS$ of
$X$, then $C\sbs \tau_2\cl\tau_1\nt C$.
\end{enumerate}
\end{proposition}

\noindent \textbf{Proof.\/} (1) Let
$$  P=Y\setminus \tau_2\cl U=Y\setminus \tau_2\cl \tau_1\nt Y=
        Y\cap \tau_2\nt\tau_1\cl(X\setminus Y). $$
Then $P\in 1\dd\Bd(X)$ and by (3) and (4) of Corollary~3.5,
$P=\tau_2\cl P=P_2^i$. Since each point of $P$ is $2$-open in $P$
and $P\in \tau_2'$, we have that $P\sbs Y_2^i$ and, hence,
$Z=Y\setminus Y_2^i\sbs \tau_2\cl U$. Clearly $Y=U\cup (Y\setminus
U)$, where
$$  \tau_1\nt(Y\setminus U)=
        \tau_1\nt \big(Y\cap (X\setminus \tau_1\nt Y)\big)=\vnth $$
and, once more applying (3) and (4) of Corollary~3.5,
$$  Y\setminus U=\tau_2\cl (Y\setminus U)=(Y\setminus U)_2^i.    $$
Thus $Y=A\cup B$, where $A=U\in\tau_1$ and $B=(Y\setminus
U)=\tau_2\cl B=B_2^i$ (clearly, if $\tau_1\nt Y=\vnth$, then
$Y=\tau_2\cl Y=Y_2^i$ and $Y=\vnth\cup Y)$.

(2) First of all, let us note that $\tau_1\nt C\neq \vnth$, since
according to (4) of Corollary~3.5 the contrary means that
$C=C_2^i$, that is, $C$ is not $2$-connected. If $C\setminus
\tau_2\cl \tau_1\nt C\neq\vnth$, then
$$  C=(C\cap \tau_2\cl\tau_1\nt C)\cup
        (C\setminus \tau_2\cl\tau_1\nt C), $$
where $A=C\cap\tau_2\cl\tau_1\nt C\in\co\tau_2'\setminus
\{\vnth\}$ in $(C,\tau_1',\tau_2')$. It is clear that  \linebreak
$C\setminus \tau_2\cl\tau_1\nt C\!\in\!1\dd\Bd(X)$ and by (3) of
Corollary~3.5, the set $B\!=\!C\setminus \tau_2\cl\tau_1\nt C\in
 \linebreak      \co\tau_2'\setminus \{\vnth\}$. Hence $C=A\cup B$, where $A,\,B\in
\co\tau_2'\setminus \{\vnth\}$ and $A\cap B=\vnth$. A
contradiction with $C$ is $2$-connected.~$\Box$ \vskip+0.2cm

Take place the following principal

\begin{theorem}{3.13}
Every $1\dd\T_1$, $p\,$-normal, $p\,$-connected and $(2,1)$-nodec
$\BS$         \linebreak           $(X,\tau_1<_C\tau_2)$ for which
every $2$-closed subset is $(1,2)\dd\WS$-supernormal in $X$, is a
$(1,2)\dd\BrS$ and, hence, a $1\dd\BrS$.
\end{theorem}

To prove this theorem, we have to formulate

\begin{lemma}{3.14}
For a $\BS$ $(X,\tau_1<\tau_2)$ the following conditions are
satisfied:
\begin{enumerate}
\item[$(1)$] If $F\in\co\tau_1$, $\Phi\in\co\tau_2$, $\Phi\sbs F$
and the $\BsS$ $(\Phi,\tau_1'',\tau_2'')$ is
$p\,\dd\WS$-supernormal in $X$, then $(\Phi,\tau_1'',\tau_2'')$ is
$p\,\dd\WS$-supernormal in $F$.

\item[$(2)$] If $F\in\co\tau_1$, $F=\tau_1\cl U$ for some $U\in
\tau_1\setminus \{\vnth\}$ and $p\,\dd\ind F=0$, then
$(X,\tau_1,\tau_2)$ is $p\,$-disconnected.
\end{enumerate}
\end{lemma}

\noindent \textbf{Proof.\/} (1) First, let $A\in \co\tau_1''$,
$B\in\co\tau_2'$ in $(F,\tau_1',\tau_2')$ and $A\cap B=\vnth$.
Clearly $B\in\co\tau_2$ and since $(\Phi,\tau_1'',\tau_2'')$ is
$(1,2)\dd\WS$-supernormal in $X$, there are sets $U\in \tau_2''$,
$V'\in \tau_1$ such that $A\sbs U$, $B\sbs V'$ and $U\cap
V'=\vnth$. Since $B\sbs F$, we have $B\sbs F\cap V'=V$, where
$V\in \tau_1'$. Therefore, $A\sbs U$, $B\sbs V$, where $U\in
\tau_2''$, $V\in \tau_1'$ and $U\cap V=\vnth$ so that
$(\Phi,\tau_1'',\tau_2'')$ is $(1,2)\dd\WS$-supernormal in $F$.

Now, let $A\in \co\tau_2''$, $B\in co\tau_1'$ and $A\cap B=\vnth$.
Because $B\in \co\tau_1$ and $(\Phi,\tau_1'',\tau_2'')$ is
$(2,1)\dd\WS$-supernormal in $X$, there are sets $U\in \tau_1''$,
$V'\in \tau_2$ such that $A\sbs U$, $B\sbs V'$ and $U\cap
V'=\vnth$. But $B\sbs F$ implies that $B\sbs F\cap V'=V\in\tau_2'$
so that $A\sbs U$, $B\sbs V$, where $U\in \tau_1''$, $V\in
\tau_2'$ and $U\cap V=\vnth$. Thus $(\Phi,\tau_1'',\tau_2'')$ is
$(2,1)\dd\WS$-supernormal in $F$.

(2) $p\,\dd\ind F=0\llra ((1,2)\dd\ind F=0$ and $(2,1)\dd\ind
F=0)$. Evidently, $(1,2)\dd\ind_xF=0$ for each point $x\in U$,
where $\tau_1\cl U=F$ and $U\in \tau_1\setminus \{\vnth\}$. Let
$U''(x)\in \tau_1'$ be any neighborhood of an arbitrarily fixed
point $x\in U$ in the $\BsS$ $(F,\tau_1',\tau_2')$. Then there is
a neighborhood $U'(x)\in \tau_1$ such that $U'(x)\cap F=U''(x)$.
If $U(x)=U'(x)\cap U$, then $U(x)\in\tau_1$ and since $U(x)\sbs
F$, we have $U(x)\in\tau_1'$. Since $(1,2)\dd\ind_xF=0$, by (1) of
Corollary~3.1.6 in [15], there is a neighborhood $V(x)\in \tau_1'$
such that $V(x)\sbs U(x)$ and $\tau_2'\cl V(x)\setminus
V(x)=\vnth$. But for $V(x)\in \tau_1'$ there is $V'(x)\in\tau_1$
such that $V'(x)\cap F=V(x)$. Because $U\sbs F$ and $V(x)\sbs
U(x)\sbs U$, we have $V(x)=V'(x)\cap U$ and so $V(x)\in\tau_1$.
Moreover, $F\in\co\tau_1\sbs\co\tau_2$ and so $V(x)=\tau_2'\cl
V(x)=\tau_2\cl V(x)$. Therefore, $V(x)\in\tau_1\cap\co\tau_2$,
$\vnth\neq V(x)\neq X$ and by (c) of Theorem~A in [18],
$(X,\tau_1,\tau_2)$ is $p\,$-disconnected.~$\Box$ \vskip+0.2cm

Now, we can proceed to prove Theorem~3.13.

Following Theorem~3.9 and Proposition~3.11, it suffices to prove
only that if $U\in\tau_1\setminus\{\vnth\}$, then $\tau_1\cl U$ is
not $(1,2)$-strongly $\sg$-discrete. Contrary: there is a set
$U\in\tau_1\setminus\{\vnth\}$ such that $\tau_1\cl U=F$ is
$(1,2)$-strongly $\sg$-discrete. Then
$F=\bigcup\limits_{k=1}^\infty F_k$, where $F_k=\tau_2\cl
F_k=(F_k)_1^i$ and since $\tau_1\sbs\tau_2$, we have
$F_k=(F_k)_2^i$ so that $p\,\dd\Ind F_k=0$ for each
$k=\ol{1,\infty}$. Furthermore, $F_k\in\co\tau_2\sbs p\,\dd\Cl(X)$
and $F_k\sbs F$ imply that $F_k\in\co\tau_2'\sbs p\,\dd\Cl(F)$ in
$(F,\tau_1',\tau_2')$ for each $k=\ol{1,\infty}$. Since
$(X,\tau_1,\tau_2)$ is $p\,$-normal, by Corollary~3.2.6 in [15],
$(F,\tau_1',\tau_2')$ is also $p\,$-normal as $F\in\co\tau_1\sbs
p\,\dd\Cl(X)$. Moreover, since each $F_k$ is $2$-closed in $X$, by
the remark between Definition~2.16 and Example~2.17, each $F_k$ is
$(2,1)\dd\WS$-supernormal in $X$ and so, by the assumption of this
theorem and (1) of Lemma~3.14, each $F_k$ is      \linebreak
$p\,\dd\WS$-supernormal in $Y$. Therefore, Corollary~2.25 gives
that $p\,\dd\Ind F=0$ and, hence, by the second part of
Corollary~3.2.8 in [15], $p\,\dd\ind F=0$ as $X$ is $1\dd\T_1$.
Thus, by (2) of Lemma~3.14, $(X,\tau_1,\tau_2)$ is
$p\,$-disconnected, which contradicts the assumption.~$\Box$

\begin{corollary}{3.15}
Every $1\dd\T_1$, $p\,$-normal, $p\,$-connected and $1$-submaximal
$($or         \linebreak        $(1,2)$-submaximal$)$ $\BS$
$(X,\tau_1<_C\tau_2)$ for which every $2$-closed subset is
\linebreak     $(1,2)\dd\WS$-supernormal in $X$, is a
$(1,2)\dd\BrS$ and, hence, a $1\dd\BrS$.
\end{corollary}

\noindent \textbf{Proof.\/} Follows directly from (1) (or (2)) of
Corollary~3.7 and Theorem~3.13.~$\Box$

\begin{definition}{3.16}
\rm A $\BS$ $(X,\tau_1,\tau_2)$ is an $(i,j)\dd\,I$-space if its
$j$-derived set is $j$-clo\-sed and $i$-discrete, that is, if
$X_j^d=\tau_j\cl X_j^d=(X_j^d)_i^i$.
\end{definition}

\begin{proposition}{3.17}
For a $1\dd\T_1$, $\BS$ $(X,\tau_1<\tau_2)$ the following
implications hold:
$$  \begin{matrix}
        (X,\tau_1,\tau_2)\;\mbox{is a $1\dd\,I$-space} & \lra &
            (X,\tau_1,\tau_2)\;\mbox{is a $(2,1)\dd\,I$-space} \\
        \big\Downarrow & & \big\Downarrow \\
        (X,\tau_1,\tau_2)\;\mbox{is a $(1,2)\dd\,I$-space} & \lra &
            (X,\tau_1,\tau_2)\;\mbox{is a $2\dd\,I$-space}.
    \end{matrix}        $$
\end{proposition}

\noindent \textbf{Proof.\/} Evidently, $X_1^d=\tau_1\cl
X_1^d=(X_1^d)_1^i$ implies that $X_1^d=\tau_1\cl
X_1^d=(X_1^d)_2^i$ and $X_2^d=\tau_2\cl X_2^d=(X_2^d)_1^i$ implies
that $X_2^d=\tau_2\cl X_2^d=(X_2^d)_2^i$ so that the horizontal
implications hold. Moreover, $X_1^d=(X_1^d)_1^i$ implies that
$X_2^d=(X_2^d)_1^i$ and $X_1^d=(X_1^d)_2^i$ implies that
$X_2^d=(X_2^d)_2^i$. But $(X,\tau_1,\tau_2)$ is also $2\dd\T_1$ so
that $X_2^d=\tau_2\cl X_2^d$ and we obtain the vertical
implications.~$\Box$

\begin{remark}{3.18}
\rm If a nonempty $\BS$ $(X,\tau_1,\tau_2)$ is an $i\dd\,I$-space,
then \linebreak     $X\,\ol{\in}\,i\dd\DI(X)$, since by [5,
p.~221], $X\in i\dd\,\ST(X)$.
\end{remark}

Hence, if a nonempty $\BS$ $(X,\tau_1<\tau_2)$ is an
$i\dd\,I$-space or an $(i,j)\dd\,I$-space, then
$X\,\ol{\in}\,2\dd\DI(X)$. But, if a $\BS$ $(X,\tau_1<\tau_2)$ is
a $(2,1)\dd\,I$-space and $X\in 1\dd\DI(X)\setminus
2\dd\DI(X)\;(\llra X\in p\,\dd\DI(X)\setminus 2\dd\DI(X)$ by (2)
of Proposition~1.4.2 in [6]), then $X=X_2^i$. Indeed: $X\in
1\dd\DI(X)\setminus 2\dd\DI(X)$ implies that $X_2^d\neq X=X_1^d$
and since $(X,\tau_1,\tau_2)$ is a $(2,1)\dd\,I$-space, we have:
$X=X_1^d=(X_1^d)_2^i=X_2^i$.

According to Proposition~3.17 above and Proposition~1.5 in [5] it
is also evident that every $1\dd\T_1$ $\BS$ $(X,\tau_1<\tau_2)$
with only finitely many $1$-nonisolated points is an
$i\dd\,I$-space and an $(i,j)\dd\,I$-space.

Proposition~3.19 as well as (2) of Proposition~3.20 and
Corollaries~3.21, 3.22 describes the connections between
topological and bitopological versions of submaximal, nodec and
$I$-spaces.

\begin{proposition}{3.19}
In the class of $1$-scattered $(\llra$ $p\,$-scattered$)$ $\BS$'s
of the type $(X,\tau_1<\tau_2)$ the following conditions are
equivalent:
\begin{enumerate}
\item[$(1)$] $(X,\tau_1,\tau_2)$ is a $1\dd\,I$-space.

\item[$(2)$] $(X,\tau_1,\tau_2)$ is $d$-submaximal and
$p\,$-submaximal.

\item[$(3)$] $(X,\tau_1,\tau_2)$ is $1$-nodec.
\end{enumerate}
\end{proposition}

\noindent \textbf{Proof.\/} In the class of $1$-scattered $\BS$'s
of the type $(X,\tau_1<\tau_2)$ the equivalences:
$(X,\tau_1,\tau_2)$ is a $1\dd\,I$-space $\llra (X,\tau_1,\tau_2)$
is $1$-submaximal $\llra (X,\tau_1,\tau_2)$ is $1$-nodec, are
given by Corollary~1.8 in [5]. The rest follows directly from the
implications after Definition~3.1.~$\Box$

\begin{proposition}{3.20}
The following conditions are satisfied:
\begin{enumerate}
\item[$(1)$] A $j\dd\T_1$ $\BS$ $(X,\tau_1,\tau_2)$ is an
$(i,j)\dd\,I$-space if and only if $A_j^d=(A_j^d)_i^i$ for each
subset $A\sbsq X$.

\item[$(2)$] If $(X,\tau_1,\tau_2)$ is $(i,j)$-nodec and $X_j^i\in
i\dd \,\cD(X)$, then $(X,\tau_1,\tau_2)$ is an
$(i,j)\dd\,I$-space.
\end{enumerate}
\end{proposition}

\noindent \textbf{Proof.\/} (1) The implication from right to left
is obvious. Let $X$ be an $(i,j)\dd\,I$-space and let $A\sbs X$ be
any subset. Then it remains to prove only that $A_j^d\sbs
(A_j^d)_i^i$. If $x\in A_j^d$, then $x\in X_j^d=(X_j^d)_i^i$ and,
hence, there is a neighborhood $U(x)\in \tau_i$ such that
$U(x)\cap X_j^d=\{x\}$. Since $x\in A_j^d$, it is clear that
$U(x)\cap A_j^d=\{x\}$ and so $x\in (A_j^d)_i^i$.

(2) Clearly, $\tau_i\cl X_j^i=X$ implies that $\tau_i\nt
X_j^d=\vnth$ and so $X_j^d\in (i,j)\dd\ND(X)$, since
$X_j^i\in\tau_j$. Thus $X_j^d=\tau_j\cl X_j^d=(X_j^d)_i^i$ as
$(X,\tau_1,\tau_2)$ is $(i,j)$-nodec.~$\Box$

\begin{corollary}{3.21}
Every $1\dd\T_1$ and $(2,1)\dd\,I$-space $(X,\tau_1<\tau_2)$ is
$(2,1)$-sub\-ma\-xi\-mal and, hence, $2$-nodec.
\end{corollary}

\noindent \textbf{Proof.\/} Let $A\sbs X$ be any set. Then, by (1)
of Proposition~3.20,
$$  \tau_1\cl A\setminus A\sbs A_1^d=(A_1^d)_2^i        $$
and so $\tau_1\cl A\setminus A=(\tau_1\cl A\setminus A)_2^i$.
Hence, by (6) of Corollary~3.5, $(X,\tau_1,\tau_2)$ is
$(2,1)$-submaximal and by the implications after Definition~3.1,
$(X,\tau_1,\tau_2)$ is $2$-submaximal. Thus, it remains to use
Corollary~1.3 in [5].~$\Box$

\begin{corollary}{3.22}
Every $1\dd\T_1$ and $(1,2)$-submaximal $\BS$ $(X,\tau_1<\tau_2)$,
for which $X_1^i\in 2\dd\,\cD(X)$, is a $(2,1)\dd\,I$-space, a
$2\dd\,I$-space and a $2$-nodec space.
\end{corollary}

\noindent \textbf{Proof.\/} Indeed, by (2) of Corollary~3.7,
$(X,\tau_1,\tau_2)$ is $(2,1)$-nodec. Hence, it remains to use (2)
of Proposition~3.20, Proposition~3.17 and Corollary~1.3 in [5],
since $(X,\tau_1,\tau_2)$ is $2$-submaximal.~$\Box$ \vskip+0.2cm

Taking into account (2) of Proposition~3.20 and Theorem~1.6 in
[5], Proposition~3.17 implies that if a $1\dd\T_1$ $\BS$
$(X,\tau_1<\tau_2)$ is an $i\dd\,I$-space or an
$(i,j)\dd\,I$-space, then $X_2^i\in 2\dd\,\cD(X)\sbs
1\dd\,\cD(X)$. Moreover, as we mentioned above (see Remark~3.18),
under the same hypotheses, $X\,\ol{\in}\,2\dd\DI(X)$. By
Theorem~1.6 in [5], for a $1\dd\T_1$ and $1\dd\,I$-space
$(X,\tau_1<\tau_2)$, the set $X_1^i\in 1\dd\,\cD(X)$. In this
context take place

\begin{proposition}{3.23}
If for a $1\dd\T_1$ and $(2,1)\dd\,I$-space $(X,\tau_1<\tau_2)$ we
have $X_1^i\in 1\dd\,\cD(X)\setminus 2\dd\,\cD(X)$, then
$\tau_2\nt X_1^d\sbs X_2^i$ and $X\,\ol{\in}\,1\dd\DI(X)\;(\llra
X\,\ol{\in}\,p\,\dd\DI(X))$.
\end{proposition}

\noindent \textbf{Proof.\/} First of all, note that
$X\,\ol{\in}\,1\dd\DI(X)$ is evident, and by (2) of
Proposition~1.4.2 in [15], $X\,\ol{\in}\,p\,\dd\DI(X)$.
Furthermore, $X_1^i\,\ol{\in}\,2\dd\,\cD(X)$ implies that
\linebreak      $X_1^d\,\ol{\in}\,2\dd\Bd(X)$ and since
$(X,\tau_1,\tau_2)$ is a $(2,1)\dd\,I$-space, we have $\vnth\neq
\tau_2\nt X_1^d=\tau_2\nt(X_1^d)_2^i$. Hence, if $x\in \tau_2\nt
X_1^d$ is any point, then there are neighborhoods $U(x),\,V(x)\in
\tau_2$ such that $U(x)\sbs X_1^d$ and $V(x)\cap X_1^d=\{x\}$.
Evidently,
$$  U(x)\cap V(x)=(U(x)\cap X_1^d)\cap V(x)=
        U(x)\cap (X_1^d\cap V(x))=\{x\}\in \tau_2       $$
and so $x\in X_2^i$.~$\Box$

\begin{proposition}{3.24}
Every $1\dd\T_1$ $\BS$ $(X,\tau_1<\tau_2)$ with only finitely many
\linebreak     $1$-nonisolated points is a $d\dd\,I$-space, a
$p\,\dd\,I$-space, $d$-submaximal, $p\,$-submaximal, $d$-nodec,
$p\,$-nodec, $d$-scattered and $p\,$-scattered.
\end{proposition}

\noindent \textbf{Proof.\/} Indeed, by Proposition~1.5 in [5],
$(X,\tau_1,\tau_2)$ is a $1\dd\,I$-space. Hence, by
Proposition~3.17, $(X,\tau_1,\tau_2)$ is a $d\dd\,I$-space and a
$p\,\dd\,I$-space. Further, by Corollary~1.3 in [5],
$(X,\tau_1,\tau_2)$ is $1$-submaximal and, hence, by the
implications after Definition~3.1, it is $d$-submaximal and
$p\,$-submaximal. Since $(X,\tau_1,\tau_2)$ is \linebreak
     $1$-submaximal, by (1) of Corollary~3.7, it is $d$-nodec and
$p\,$-nodec. Finally, by [5, p.~221], $(X,\tau_1,\tau_2)$ is
$d$-scattered, since $(X,\tau_1,\tau_2)$ is a $d\dd\,I$-space and,
therefore, by (3) of Proposition~1.4.12 in [15],
$(X,\tau_1,\tau_2)$ is $p\,$-scattered.~$\Box$

\begin{proposition}{3.25}
For a $\BS$ $(X,\tau_1<\tau_2)$ the following conditions are
satisfied:
\begin{enumerate}
\item[$(1)$] If $(X,\tau_1<_C\tau_2)$ is $(1,2)$-submaximal, then
$X=Y\cup Z$, where $Y\in \co\tau_2$, $(Y,\tau_1',\tau_2')$ is a
$(2,1)\dd\,I$-space and $Z\in \co\tau_1\cap 1\dd\DI(Z)$,
$(Z,\tau_1'',\tau_2'')$ is   \linebreak       $(1,2)$-submaximal.

\item[$(2)$] If $(X,\tau_1,\tau_2)$ is $(1,2)$-nodec, then
$X=Y\cup Z$, where $Y\in \co\tau_1$, $(Y,\tau_1',\tau_2')$ is a
$(1,2)\dd\,I$-space and $Z\in \co\tau_2\cap 2\dd\DI(Z)$,
$(Z,\tau_1'',\tau_2'')$ is $(1,2)$-nodec.
\end{enumerate}
\end{proposition}

\noindent \textbf{Proof.\/} (1) It is clear that $X_1^i\in
2\dd\,\cD(\tau_2\cl X_1^i)$. If $Y=\tau_2\cl X_1^i$, then by
Corollary~3.6, $(Y,\tau_1',\tau_2')$ is also $(1,2)$-submaximal
and, hence, by (2) of Corollary~3.7, $(Y,\tau_1',\tau_2')$ is
$(2,1)$-nodec. On the other hand, $Y_1^i\in 2\dd\,\cD(Y)$ and
according to (2) of Proposition~3.20, $(Y,\tau_1',\tau_2')$ is a
$(2,1)\dd\,I$-space. Furthermore, since $\tau_1<_C\tau_2$ and
$X_1^i\in \tau_1$, by (3) of Corollary~2.2.7 in [15], $Y=\tau_2\cl
X_1^i=\tau_1\cl X_1^i$. Therefore, $X\setminus Y=X\setminus
\tau_1\cl X_1^i\in \tau_1$ and, hence, $(X\setminus Y)\cap
X_1^i=\vnth$ implies that $(X\setminus Y)_1^i=\vnth$, so that
$(\tau_1\cl (X\setminus Y))_1^i=\vnth$.

Clearly $X=Y\cup Z$, where $Z=\tau_1\cl (X\setminus Y)$ and once
more applying Corollary~3.6, we obtain that
$(Z,\tau_1'',\tau_2'')$ is also $(1,2)$-submaximal,
$Z=\co\tau_1\cap 1\dd\DI(Z)$.

(2) Evidently, $X_2^i\in 1\dd\cD(\tau_1\cl X_2^i)$. If
$Y=\tau_1\cl X_2^i$, then by Remark~1.2, $(Y,\tau_1',\tau_2')$ is
also $(1,2)$-nodec. Hence, according to (2) of Proposition~3.20,
$(Y,\tau_1',\tau_2')$ is $(1,2)\dd\,I$-space, since $Y_2^i\in
1\dd\,\cD(Y)$. Clearly $X\setminus Y\in \tau_1\sbs \tau_2$ and
since $(X\setminus Y)\cap X_2^i=\vnth$, we have $(X\setminus
Y)_2^i=\vnth$.

Therefore, $(\tau_2\cl(X\setminus Y))_2^i=\vnth$ and if
$Z=\tau_2\cl(X\setminus Y)$, then $Z\in \co\tau_2\cap 2\dd\DI(Z)$
and, by Remark 1.2, $(Z,\tau_1'',\tau_2'')$ is $(1,2)$-nodec.
$\Box$

\begin{definition}{3.26}
\rm A surjection $f: (X,\tau_1,\tau_2)\to (Y,\gm_1,\gm_2)$ is
$(i,j)$-locally closed quotient (briefly, $(i,j)\dd\lcq)$, if
$f^{-1}(A)\in (i,j)\dd\LC(X)$ implies that $A\in (i,j)\dd\LC(Y)$
or, equivalently, if $f^{-1}(A)\in\co(i,j)\dd\LC(X)$ implies that
$A\in \co(i,j)\dd\LC(Y)$ for each subset $A\sbs Y$.
\end{definition}

\begin{proposition}{3.27}
For a function $f: (X,\tau_1,\tau_2)\to (Y,\gm_1,\gm_2)$ the
following conditions are satisfied:
\begin{enumerate}
\item[$(1)$] If $f$ is a $d$-open or a $d$-closed surjection, then
$f$ is $(i,j)\dd\lcq$.

\item[$(2)$] If $f$ is $(i,j)\dd\lcq$ and $(X,\tau_1,\tau_2)$ is
$(i,j)$-submaximal, then $(Y,\gm_1,\gm_2)$ is also
$(i,j)$-submaximal.
\end{enumerate}
\end{proposition}

\noindent \textbf{Proof.\/} (1) It suffices to consider only the
case where $f$ is a $d$-open surjection. Let $B\sbs Y$ and
$f^{-1}(B)\in \co(i,j)\dd\LC(X)$, that is, $f^{-1}(B)=A_1\cup
A_2$, where $A_1\in \tau_j$ and $A_2\in\co\tau_i$. Then
$B_1=f(A_1)\in\gm_j$, $B_2=Y\setminus f(X\setminus
A_2)\in\co\gm_i$ and $B=B_1\cup B_2$.

(2) Let $A\sbs Y$ be any subset. Since $(X,\tau_1,\tau_2)$ is
$(i,j)$-submaximal, $f^{-1}(A)\in (i,j)\dd\LC(X)$ and by
Definition~3.26, $A\in (i,j)\dd\LC(Y)$ so that $(Y,\gm_1,\gm_2)$
is        \linebreak       $(i,j)$-submaximal.~$\Box$

\begin{corollary}{3.28}
If $f: (X,\tau_1,\tau_2)\to (Y,\gm_1,\gm_2)$ is a $d$-open or a
$d$-closed surjection and $(X,\tau_1,\tau_2)$ is
$p\,$-submaximal, then $(Y,\gm_1,\gm_2)$ is also $p\,$-submaximal.
\end{corollary}

\begin{corollary}{3.29}
If $(X,\tau_1,\tau_2)$ is $p\,$-submaximal and $\tau_i\sbsq\gm_i$,
then $(X,\gm_1,\gm_2)$ is also $p\,$-submaximal.
\end{corollary}

\noindent \textbf{Proof.\/} Indeed, the identity map $i_x:
(X,\tau_1,\tau_2)\to (X,\gm_1,\gm_2)$ is a $d$-open and $d$-closed
surjection.~$\Box$

\begin{theorem}{3.30}
If a $\BS$ $(X,\tau_1<\tau_2)$ is a $(2,1)\dd\,I$-space and $f:
(X,\tau_1<\tau_2)\to (Y,\gm_1<\gm_2)$ is a $d$-closed or a
$d$-open surjection, then $(Y,\gm_1<\gm_2)$ is also a
$(2,1)\dd\,I$-space.
\end{theorem}

\noindent \textbf{Proof.\/} First, let $f$ be a $d$-closed
surjection. Since $(X,\tau_1,\tau_2)$ is a $(2,1)\dd\,I$-space,
for the set $F=X_1^d=\tau_1\cl X_1^d=(X_1^d)_2^i$ the $\BsS$
$(F,\tau_1'<\tau_2')$ is $2$-discrete. Hence, $2^F\sbs
\co\tau_2'\sbs\co\tau_2$ and, also, $\co\tau_1'\sbs\co\tau_1$ as
$F\in \co\tau_1\sbs\co\tau_2$. Therefore, the restriction
$f|_{{}_F}: (F,\tau_1',\tau_2')\to (Y,\gm_1,\gm_2)$ is $d$-closed.
Let $P=f(F)\sbs Y$. Evidently, for each subset $B\sbsq P$ there is
a subset $A\sbsq F$ such that $f(A)=B$ so that each subset $B\sbsq
P$ is $2$-closed in $(Y,\gm_1,\gm_2)$ as each subset $A\sbsq F$ is
$2$-closed in $(X,\tau_1,\tau_2)$. Hence, the $\BsS$
$(P,\gm_1'<\gm_2')$ is also $2$-discrete so that $P=\gm_1\cl
P=P_2^i$ as $P=f(F)$, $F\in \co\tau_1$ and $f$ is $1$-closed.
Therefore, to completes the proof of the first part, it suffices
to prove only that each point $x\in Y\setminus P$ is $1$-isolated
in $(Y,\gm_1,\gm_2)$.

Let $x\in Y\setminus P$ be any point. Then $E=f^{-1}(x)$ implies
that $E\sbs X\setminus F=X\setminus X_1^d=X_1^i$. Hence,
$E\in\tau_1$ and so $\Phi=X\setminus E\in\co\tau_1$. Because $f$
is $1$-closed, the set $f(\Phi)\in\co\gm_1$ and $\{x\}=Y\setminus
f(\Phi)$, that is, $\{x\}\in \gm_1$ so that $x\in Y_1^i$.
Therefore, $Y_1^d=P$ and, hence, $Y_1^d=\tau_1\cl
Y_1^d=(Y_1^d)_2^i$, that is, $(Y,\gm_1,\gm_2)$ is a
$(2,1)\dd\,I$-space.

For a $d$-open $f$, if $U=X_1^i$, then $U\in \tau_1$ and,
therefore, $f(U)\in\gm_1$, $f(U)\sbs Y_1^i$. Putting now
$P=Y\setminus f(U)$ we obtain that $P\in \co\gm_1$ and
$f^{-1}(P)\sbs F$, where $F=X\setminus U\in\co\tau_1$. Since
$(X,\tau_1,\tau_2)$ is a $(2,1)\dd\,I$-space, for the set
$F=X_1^d$, we have $F=\tau_1\cl F=F_2^i$. Therefore, all subsets
of $F$ are $2$-closed in $F$ and, hence, in $X$. Since
$f^{-1}(P)\sbs F$, for each subset $B\sbsq P$ we have
$f^{-1}(B)\in\co\tau_2$. Hence, for each subset $B\sbsq P$ the set
$X\setminus f^{-1}(B)\in\tau_2$ and since $f$ is $2$-open, the set
$f(X\setminus f^{-1}(B))=Y\setminus B\in \gm_2$. Therefore,
$P=\gm_1\cl P=P_2^i$, where each point of $Y\setminus P=f(U)$ is
$1$-isolated in $Y$. Thus $Y_1^d=P$ and $Y_1^d=\tau_1\cl
Y_1^d=(Y_1^d)_2^i$ so that $(Y,\gm_1,\gm_2)$ is a
$(2,1)\dd\,I$-space.~$\Box$ \vskip+0.2cm

The rest of this section is concerned with the notion of $D$-space
and its bitopological modifications. The sufficient conditions are
established under which a $\TS$ is a $D$-space on the one hand, by
the topological methods and, on the other hand, by the
bitopological ones.

\begin{definition}{3.31}
\rm A neighborhood assignment on a $\TS$ $(X,\tau)$ is a function
$\phi: X\to \tau$ such that $x\in \phi(x)$. A $\TS$ $(X,\tau)$ is
a $D$-space if for every neighborhood assignment $\phi$ on $X$
there is a closed discrete subset $D\!\sbsq\!X$ such that
$\bigcup\limits_{x\in D} \phi(x)\!=~\!\!X$~[21].
\end{definition}

\begin{remark}{3.32}
\rm In general, a $D$-space is not compact. Indeed, if $\phi: X\to
\tau$ is any neighborhood assignment on $X$, then there is a
subset $D=\ol{D}=D^i$ such that $\bigcup\limits_{x\in D}
\phi(x)=X$. Let us consider the open covering $\cU$ of the $\TS$
$(X,\tau)$:
$$  \cU=\big\{\phi(x)\setminus D:x\in X\setminus D\big\}\cup
        \big\{U(x):\;U(x)\in\tau,\;U(x)\cap D=\{x\}\big\}.      $$
Clearly, if $D$ is infinite, then $\cU$ does not contain a finite
subcovering.
\end{remark}

Let $\phi: X\to \tau$ be a neighborhood assignment on the $\TS$
$(X,\tau)$ and let          \linebreak      $f: (X,\tau)\to
(Y,\gm)$ be an open surjection. Then $\psi: Y\to \gm$, defined as
$\psi(y)=f\big(\bigcup\limits_{x\in f^{-1}(y)}\phi(x)\big)$, is
the neighborhood assignment on the $\TS$ $(Y,\gm)$. Also, if
$\psi: Y\to \gm$ is a neighborhood assignment on the $\TS$
$(Y,\gm)$ and $f: (X,\tau)\to (Y,\gm)$ is any continuous function,
then $\phi: X\to \tau$, defined as $\phi(x)=f^{-1}(\psi(f(x))$, is
the neighborhood assignment on the $\TS$ $(X,\tau)$.

\begin{definition}{3.33}
\rm Let $\phi: X\to \tau$ and $\psi: Y\to \gm$ be neighborhood
assignments on the $\TS$'s $(X,\tau)$ and $(Y,\gm)$, respectively.
Then we shall say that a surjection $f: X\to Y$ connects $\phi$
with $\psi$ if $\psi\circ f=f\circ \phi$, i.e., if
$\psi(f(x))=f(\phi(x))$ for each point $x\in X$.
\end{definition}

\begin{remark}{3.34}
\rm (a) If $\phi: X\to \tau$ and $\psi: Y\to \gm$ are neighborhood
assignments on the $\TS$'s $(X,\tau)$ and $(Y,\gm)$, respectively,
where $\psi$ is defined by an open bijection $f: (X,\tau)\to
(Y,\gm)$, then $f$ connects $\phi$ with $\psi$.

Indeed,
$$  \psi(f(x))=f\big({\textstyle \bigcup\limits_{x\in f^{-1}(f(x))}}
        \phi(x)\big)=f(\phi(x)) \;\;\mbox{for each}\;\; x\in X. $$

(b) If $\phi: X\to \tau$ and $\psi: Y\to \gm$ are neighborhood
assignments on the $\TS$'s $(X,\tau)$ and $(Y,\gm)$, respectively,
where $\psi$ is defined by a continuous surjection $f: (X,\tau)\to
(Y,\gm)$, then $f$ connects $\phi$ with $\psi$.

Indeed,
$$  f(\phi(x))=f\big(f^{-1}(\psi(f(x))\big)=\psi(f(x))
            \;\;\mbox{for each}\;\; x\in X. $$
\end{remark}

It is clear that if a surjection $f: X\to Y$ connects a
neighborhood assignment $\phi: X\to \tau$ with a neighborhood
assignment $\psi: Y\to \gm$, then $f(\phi(x))\in \gm$ for each
point $x\in X$; but, in general, $f$ is not open.

\begin{example}{3.35}
\rm Let $X=\{a,b,c,d,e\}$, $\tau$ be the discrete topology on $X$
and let a neighborhood assignment $\phi: X\to \tau$ on the $\TS$
$(X,\tau)$ be defined in the manner as follows:
$$  \phi(a)=\{a\}, \;\; \phi(b)=\{a,b\}, \;\; \phi(c)=\{a,b,c\},
            \;\; \phi(d)=\{a,b,c,d\}    $$
and $\phi(e)=X$. Furthermore, let $Y=\{0,1,2,3,4\}$,
$$  \gm=\{\vnth,\{0\},\{0,1\},\{0,1,2\},\{0,1,2,3\},Y\} $$
and let a neighborhood assignment $\psi:Y\to \gm$ on the $\TS$
$(Y,\gm)$ be defined as follows:
\begin{eqnarray*}
    &\ds \psi(0)=\{0\}, \;\; \psi(1)=\{0,1\}, \;\;
        \psi(2)=\{0,1,2\}, \\
    &\ds \psi(3)=\{0,1,2,3\} \;\;\mbox{and}\;\; \psi(4)=Y.
\end{eqnarray*}
If $f: X\to Y$ is defined as follows:
$$  f(a)=0, \;\; f(b)=1, \;\; f(c)=2, \;\; f(d)=3 \;\;\mbox{and} \;\;
        f(e)=4,     $$
then $f$ is a continuous bijection (i.e., a compression map), as
$\tau$ is discrete; moreover, $f$ connects $\phi$ with $\psi$.
But, $f$ is not open, since, for example,
$f(a,c)=\{0,2\}\,\ol{\in}\,\gm$, where $\{a,c\}\in\tau$, and so
$f$ is not a homeomorphism.

Now, if replace the sets $X$ and $Y$ and if we consider the map
$f^{-1}:Y\to X$, then $f^{-1}$ will become an open bijection,
connecting $\psi$ with $\phi$. Clearly, $f^{-1}$ is not
homeomorphism, since it is not continuous.

Moreover, it follows immediately from the determination of the
topology $\gm$ on the set $Y$ that for any neighborhood assignment
$\psi: Y\to \gm$ always $\psi(4)=Y$, since the set
$\{0,1,2,3,4\}=Y$ is a unique set from $\gm$ which contains the
point $4$. Hence, for any neighborhood assignment $\psi: Y\to \gm$
we can associate the closed discrete set $\{4\}$ such that
$\psi(4)=Y$. Thus $(Y,\gm)$ is a $D$-space.

Note also here that if a compression map $g: X\to Y$ is defined as
follows: $g(a)=1$, $g(b)=0$, $g(c)=2$, $g(d)=3$ and $g(e)=4$, then
$g$ does not connect $\phi$ with $\psi$ since
$\psi(g(a))=\psi(1)=\{0,1\}\neq g(\phi(a))=g(a)=\{1\}$.
\end{example}

\begin{theorem}{3.36}
The following conditions are satisfied:
\begin{enumerate}
\item[$(1)$] If for each neighborhood assignment $\phi: X\to \tau$
on the $\TS$ $(X,\tau)$ there are a neighborhood assignment
$\psi:Y\to \gm$ on the $\TS$ $(Y,\gm)$ and a compression map $f:
(X,\tau)\to (Y,\gm)$ which connects $\phi$ with $\psi$, then
$(Y,\gm)$ is a $D$-space implies that $(X,\tau)$ is also a
$D$-space.

\item[$(2)$] If for each neighborhood assignment $\psi: Y\to \gm$
on the $\TS$ $(Y,\gm)$ there are a neighborhood assignment
$\phi:X\to \tau$ on the $\TS$ $(X,\tau)$ and an open bijection $f:
(X,\tau)\to (Y,\gm)$ which connects $\phi$ with $\psi$, then
$(X,\tau)$ is a $D$-space implies that $(Y,\gm)$ is also a
$D$-space.
\end{enumerate}
\end{theorem}

\noindent \textbf{Proof.\/} (1) Let $\phi:X\to \tau$ be any
neighborhood assignment on the $\TS$ $(X,\tau)$. Then there are a
neighborhood assignment $\psi:Y\to \gm$ on the $\TS$ $(Y,\gm)$ and
a compression map $f: (X,\tau)\to (Y,\gm)$ such that $\psi\circ
f=f\circ \phi$. Since $(Y,\gm)$ is a $D$-space, for $\psi$ there
exists a closed discrete set $D'$ such that $\bigcup\limits_{y\in
D'} \psi(y)=Y$. Clearly $f^{-1}(D')=D=\ol{D}$ as $f$ is
continuous. On the other hand, if $x\in D$ is any point, then
$y=f(x)\in D'$ and so there is a neighborhood $U(y)\in\gm$ such
that $U(y)\cap D'=\{y\}$. Since $f$ is a compression map,
$$  x=f^{-1}(y)=f^{-1}(U(y)\cap D')=f^{-1}(U(y))\cap
            f^{-1}(D')=U(x)\cap D,      $$
where $U(x)\in\tau$. Hence, $D$ is closed discrete subset of
$(X,\tau)$. Finally, $\bigcup\limits_{y\in D'} \psi(y)=Y$ implies
that
\begin{eqnarray*}
    &\ds X=f^{-1}(Y)= \\
    &\ds =f^{-1}\big({\textstyle\bigcup\limits_{y\in D'}} \psi(y)\big)=
            {\textstyle\bigcup\limits_{x\in D}} f^{-1}\big(\psi(f(x))\big)=
        {\textstyle\bigcup\limits_{x\in D}} f^{-1}(f(\phi(x)))=
            {\textstyle\bigcup\limits_{x\in D}} \phi(x)
\end{eqnarray*}
as $f$ connects $\phi$ with $\psi$.

(2) If $\psi: Y\to \gm$ is any neighborhood assignment on the
$\TS$ $(Y,\tau)$, then there are a neighborhood assignment $\phi:
X\to \tau$ on the $\TS$ $(X,\tau)$ and an open bijection $f:
(X,\tau)\to (Y,\gm)$ which connects $\phi$ with $\psi$. Since
$(X,\tau)$ is a $D$-space, there is a closed discrete set $D$ such
that $\bigcup\limits_{x\in D} \phi(x)=X$. If $U=X\setminus D$,
then $D'=f(D)=f(X\setminus U)=Y\setminus f(U)\in\co\gm$. For any
point $y\in D'$ there is a neighborhood $U(x)\in\tau$ of the point
$x=f^{-1}(y)$ such that $U(x)\cap D=\{x\}$. Hence
\begin{eqnarray*}
    &\ds y=f(f^{-1}(y))=f(x)=f(U(x)\cap D)= \\
    & \ds =f\big(U(x)\cap f^{-1}(D')\big)=
        f(U(x))\cap D'=U(y)\cap D',
\end{eqnarray*}
where $U(y)\in \gm$ as $f$ is an open bijection. Therefore, $D'$
is a closed discrete subset of $(Y,\gm,)$. Finally,
$\bigcup\limits_{x\in D} \phi(x)=X$ implies that
$$  Y=f(X)=f\big({\textstyle\bigcup\limits_{x\in D}} \phi(x)\big)=
            {\textstyle\bigcup\limits_{x\in D}} f(\phi(x))=
        {\textstyle\bigcup\limits_{x\in D}} \psi(f(x))=
            {\textstyle\bigcup\limits_{y\in D'}} \psi(y)  $$
as $f$ connects $\phi$ with $\psi$.~$\Box$

\begin{theorem}{3.37}
Every connected $I$-space is a $D$-space.
\end{theorem}

\noindent \textbf{Proof.\/} By condition, $X^d=\ol{X^d}=(X^d)^i$.
Let us prove that ${\textstyle\bigcup\limits_{x\in D}} \phi(x)=X$
for any neighborhood assignment $\phi:X\to \tau$ on the $\TS$
$(X,\tau)$, where $D=X^d$. Contrary: there is a neighborhood
assignment $\phi: X\to \tau$ on the $\TS$ $(X,\tau)$ such that
${\textstyle\bigcup\limits_{x\in D}} \phi(x)\neq X$, i.e.,
$X\setminus {\textstyle\bigcup\limits_{x\in D}} \phi(x)\neq
\vnth$. Since the set $D=X^d\sbs {\textstyle\bigcup\limits_{x\in
D}} \phi(x)$, we have that
$$  X\setminus {\textstyle\bigcup\limits_{x\in D}}\phi(x)=
        X^i\setminus {\textstyle\bigcup\limits_{x\in D}}\phi(x).   $$
Therefore, it is not difficult to see, that
$$  X^i\setminus {\textstyle\bigcup\limits_{x\in D}}\phi(x)=
        X^i\setminus \ol{{\textstyle\bigcup\limits_{x\in D}}\phi(x)} $$
and, hence, for the set $A=X\setminus
{\textstyle\bigcup\limits_{x\in D}}\phi(x)\in\tau\cap \co\tau$ we
have: $\vnth\neq A\neq X$, which is impossible as $(X,\tau)$ is
connected.~$\Box$

\begin{corollary}{3.38}
Every connected $\TS$ $(X,\tau)$, for which $X^i\in \cD(X)$, is a
$D$-space if it is nodec or submaximal.
\end{corollary}

\noindent \textbf{Proof.\/} Follows directly from Theorem~1.6 and
Corollary~1.7 in [5].~$\Box$

\begin{definition}{3.39}
A $\BS$ $(X,\tau_1,\tau_2)$ will be called an $(i,j)\dd\,D$-space
if for every $i$-neighborhood assignment $\phi_i: X\to \tau_i$ on
the $\BS$ $(X,\tau_1,\tau_2)$ there is a $j$-closed $i$-discrete
subset $D\sbsq X$ such that $\bigcup\limits_{x\in D}\phi_i(x)=X$.
\end{definition}

It is clear that for a $\BS$ $(X,\tau_1<\tau_2)$ we have:
$$  (X,\tau_1,\tau_2)\;\mbox{is a $1\dd\,D$-space}\lra
        (X,\tau_1,\tau_2)\;\mbox{is $(1,2)\dd\,D$-space}  $$
and
$$  \ \hskip-1cm (X,\tau_1,\tau_2)\;\mbox{is a $(2,1)\dd\,D$-space}\lra
        (X,\tau_1,\tau_2)\;\mbox{is a $2\dd\,D$-space}.   $$

By analogy with Theorem~3.36 one can prove

\begin{theorem}{3.40}
The following conditions are satisfied:
\begin{enumerate}
\item[$(1)$] If for each $i$-neighborhood assignment $\phi_i: X\to
\tau_i$ on the $\BS$ $(X,\tau_1,\tau_2)$ there are
$i$-neighborhood assignment $\psi_i:Y\to \gm_i$ on the $\BS$
$(Y,\gm_1,\gm_2)$ and a $d$-compression map $f:
(X,\tau_1,\tau_2)\to (Y,\gm_1,\gm_2)$ which connects $\phi_i$ with
$\psi_i$, then $(Y,\gm_1,\gm_2)$ is an $(i,j)\dd\,D$-space implies
that $(X,\tau_1,\tau_2)$ is also an \linebreak
   $(i,j)\dd\,D$-space.

\item[$(2)$] If for each $i$-neighborhood assignment $\psi_i: Y\to
\gm_i$ on the $\BS$ $(Y,\gm_1,\gm_2)$ there are $i$-neighborhood
assignment $\phi_i:X\to \tau_i$ on the $\BS$ $(X,\tau_1,\tau_2)$
and a $d$-open bijection $f: (X,\tau_1,\tau_2)\to (Y,\gm_1,\gm_2)$
which connects $\phi_i$ with $\psi_i$, then $(X,\tau_1,\tau_2)$ is
an $(i,j)\dd\,D$-space implies that $(Y,\gm_1,\gm_2)$ is also an
\linebreak           $(i,j)\dd\,D$-space.
\end{enumerate}
\end{theorem}

\begin{corollary}{3.41}
The following conditions are satisfied:
\begin{enumerate}
\item[$(1)$] If for a pair $(\phi_1,\phi_2)$, where $\phi_i: X\to
\tau_i$ is an $i$-neighborhood assignment on the $\BS$
$(X,\tau_1,\tau_2)$, there are a pair $(\psi_1,\psi_2)$, where
$\psi_i:Y\to \gm_i$ is an $i$-neighborhood assignment on the $\BS$
$(Y,\gm_1,\gm_2)$, and a $d$-compression map $f:
(X,\tau_1,\tau_2)\to (Y,\gm_1,\gm_2)$ which connects
$(\phi_1,\phi_2)$ with $(\psi_1,\psi_2)$ $($i.e., $f$ connects
$\phi_1$ with $\psi_1$ and $\phi_2$ with $\psi_2)$, then
$(Y,\gm_1,\gm_2)$ is a $p\,\dd\,D$-space implies that
$(X,\tau_1,\tau_2)$ is also a $p\,\dd\,D$-space.

\item[$(2)$] If for a pair $(\psi_1,\psi_2)$, where $\psi_i: Y\to
\gm_i$ is an $i$-neighborhood assignment on the $\BS$
$(Y,\gm_1,\gm_2)$, there are a pair $(\phi_1,\phi_2)$, where
$\phi_i: X\to \tau_i$ is an $i$-neighborhood assignment on the
$\BS$ $(X,\tau_1,\tau_2)$, and a $d$-open bijection \linebreak
  $f: (X,\tau_1,\tau_2)\to (Y,\gm_1,\gm_2)$ which connects
$(\phi_1,\phi_2)$ with $(\psi_1,\psi_2)$, then $(X,\tau_1,\tau_2)$
is a $p\,\dd\,D$-space implies that $(Y,\gm_1,\gm_2)$ is also a
$p\,\dd\,D$-space.
\end{enumerate}
\end{corollary}

\begin{theorem}{3.42}
Every $p\,$-connected $(i,j)\dd\cI$-space is an
$(i,j)\dd\,D$-space.
\end{theorem}

\noindent \textbf{Proof.\/} By condition, $X_j^d=\tau_j\cl
X_j^d=(X_j^d)_i^i$. Let us prove that $\bigcup\limits_{x\in D}
\phi(x)=X$ for any $i$-neighborhood assignment $\phi_i:X\to
\tau_i$ on the $\BS$ $(X,\tau_1,\tau_2)$, where $D=X_j^d$. Indeed,
if there is an $i$-neighborhood assignment $\phi_i:X\to \tau_i$ on
the $\BS$ $(X,\tau_1,\tau_2)$ such that $\bigcup\limits_{x\in D}
\phi_i(x)\neq X$, then
$$  \vnth\neq X\setminus {\textstyle\bigcup\limits_{x\in D}} \phi_i(x)=
        X\setminus \tau_j\cl{\textstyle\bigcup\limits_{x\in D}}\phi_i(x)  $$
as
$$  X\setminus {\textstyle\bigcup\limits_{x\in D}}\phi_i(x)=
        X_j^i\setminus {\textstyle\bigcup\limits_{x\in D}}\phi_i(x).   $$
Therefore, if $A=X\setminus {\textstyle\bigcup\limits_{x\in
D}}\phi_i(x)$, then $A\in \tau_j\cap \co\tau_i$, $\vnth\neq A\neq
X$ and so, by (c) of Theorem~A in [18], $(X,\tau_1,\tau_2)$ is not
$p\,$-connected.~$\Box$

\begin{corollary}{3.43}
Every $p\,$-connected $p\dd\cI$-space is a $p\,\dd\,D$-space.
\end{corollary}

\begin{corollary}{3.44}
Every $p\,$-connected $(i,j)$-nodec $\BS$ $(X,\tau_1,\tau_2)$, for
which $X_j^i\in i\dd\,\cD(X)$, is an $(i,j)\dd\,D$-space.
\end{corollary}

\noindent \textbf{Proof.\/} Follows directly from (2) of
Proposition~3.20.~$\Box$ \vskip+0.2cm

Let $\Phi=\{\phi\}$ be the family of all neighborhood assignments
on a $\TS$ $(X,\tau)$. Then a binary relation $\leq$, defined on
$\Phi$ in the manner as follows: $\phi_1,\,\phi_2\in \Phi$,
$\phi_1\leq \phi_2$ if $\phi_1(x)\sbsq \phi_2(x)$ for each point
$x\in X$, is a partial order on $\Phi$. Evidently, this partial
order is linear if for each pair $\phi_1,\,\phi_2\in \Phi$ and
each point $x\in X$ we have $\phi_1(x)\sbsq \phi_2(x)$ or
$\phi_2(x)\sbsq \phi_1(x)$.

\begin{theorem}{3.45}
If for a $\TS$ $(X,\tau)$ the family $\Phi$ of all neighborhood
assignments $\phi:X\to \tau$ on the $\TS$ $(X,\tau)$ is linearly
ordered and there are topologies $\tau_1$ and $\tau_2$ on $X$ such
that $\sup(\tau_1,\tau_2)=\tau$ and $(X,\tau_1,\tau_2)$ is a
$p\,\dd\,D$-space, then $(X,\tau)$ is a $D$-space.
\end{theorem}

\noindent \textbf{Proof.\/} Let $\phi:X\to \tau$ be any
neighborhood assignment on the $\TS$ $(X,\tau)$, where
$\tau=\sup(\tau_1,\tau_2)$ for some topologies $\tau_1$ and
$\tau_2$ on $X$. Without loss of generality we can suppose that
for each point $x\in X$ the set $\phi(x)$ is basic open. Then
$\phi(x)=U_1(x)\cap U_2(x)$, where $x\in U_i(x)\in\tau_i$. Let us
define the neighborhood assignments $\phi_i:X\to\tau_i$ as
follows: $\phi_i(x)=U_i(x)$ for each point $x\in X$. Since
$(X,\tau_1,\tau_2)$ is a $p\,\dd\,D$-space, for $\phi_1$ there
exists a $2$-closed $1$-discrete set $D_1$ and for $\phi_2$ there
exists a $1$-closed $2$-discrete set $D_2$ such that
$$  \bigcup\limits_{x\in D_1} \phi_1(x)=X=
        \bigcup\limits_{x\in D_2} \phi_2(x).        $$
It is clear, that $\phi_1,\,\phi_2\in \Phi$ as $\tau_1\cup
\tau_2\sbs \tau$.

Furthermore, let $D=D_1\cup D_2$. Because $\tau_1\cup \tau_2\sbs
\tau$, we have $D_1\in\co\tau_2\sbs\co\tau$,
$D_2\in\co\tau_1\sbs\co\tau$ and, hence, $D\in\co\tau$. Moreover,
let $x\in D$ be any point, where
$$  D=(D_1\setminus D_2)\cup (D_1\cap D_2)\cup (D_2\setminus D_1).  $$
If $x\in D_1\cap D_2$, then $D_1=(D_1)_1^i$ and $D_2=(D_2)_2^i$
imply that there are neighborhoods $U(x)\in\tau_1$ and
$V(x)\in\tau_2$ such that $U(x)\cap D_1=\{x\}=V(x)\cap D_2$. Since
$U(x)\in \tau_1\sbs \tau$, $V(x)\in \tau_2\sbs \tau$, we have that
$W(x)=U(x)\cap V(x)\in \tau$ and $W(x)\cap D=\{x\}$ so that
$D_1\cap D_2\sbs D^i$. Now, if $x\in D_i\setminus D_j$, then there
is a neighborhood $U(x)\in \tau_i\sbs \tau$ such that $U(x)\cap
D_i=\{x\}$. Let $V(x)=U(x)\setminus D_j$. Then $V(x)\in \tau_i\sbs
\tau$ and $V(x)\cap D=\{x\}$, so that $D_i\setminus D_j\sbs D^i$
too.

Therefore, $D=D^i$ and it remains to prove only that
$\bigcup\limits_{x\in D} \phi(x)=X$. But
$$  {\textstyle \bigcup\limits_{x\in D}} \phi(x)=
        {\textstyle \bigcup\limits_{x\in D}}
            \big(\phi_1(x)\cap \phi_2(x)\big).     $$
Since $\phi_1,\,\phi_2\in \Phi$ and $\Phi$ is linearly ordered, we
have that $\phi_1(x)\sbsq \phi_2(x)$ or $\phi_2(x)\sbsq \phi_1(x)$
for each point $x\in D$. Let, for example, $\phi_2(x)\sbsq
\phi_1(x)$ for each point $x\in D$. Then
$$  {\textstyle \bigcup\limits_{x\in D}} \phi(x)=
        {\textstyle \bigcup\limits_{x\in D}} \phi_2(x)=X    $$
and, consequently, $(X,\tau)$ is a $D$-space.~$\Box$

\vskip+0.8cm

\end{document}